\documentclass[final,onefignum,onetabnum,onethmnum]{siamltex}
\pdfoutput=1

\usepackage{multicol}
\usepackage{amsfonts}
\usepackage{graphicx}

\newcommand{\BC}{\mathbb{D}}
\newcommand{\HC}{\mathbb{H}}
\newcommand{\C}{\mathbb{C}}
\newcommand{\N}{\mathbb{N}}

\font\ff=eusm10 scaled 1200
\newcommand{\K}{\hbox{\ff K}}
\newcommand{\Cap}{\mathop{\mathrm{cap}}}
\newcommand{\re}{\mathop{\mathrm{Re}}}
\newcommand{\im}{\mathop{\mathrm{Im}}}
\newcommand{\M}{\mbox{\sf M}}

\def\QM{\mathrm{QM}}

\graphicspath{{PDF/}{EPS/}}  

\usepackage{amsmath}

\title{On moduli of rings and quadrilaterals: algorithms and experiments}

\author{Harri Hakula\thanks{Aalto University, Institute of Mathematics,
        P.O. Box 11100, FI-00076 Aalto, Finland ({\tt harri.hakula@tkk.fi})}
       \and Antti Rasila\thanks{Aalto University, Institute of Mathematics,
        P.O. Box 11100, FI-00076 Aalto, Finland ({\tt antti.rasila@iki.fi})}
       \and Matti Vuorinen\thanks{Department of Mathematics,
        FI-20014 University of Turku, Finland ({\tt vuorinen@utu.fi})}}

\begin{document}

\maketitle

\begin{abstract}
Moduli of rings and quadrilaterals are frequently
applied in geometric function theory, see e.g. the Handbook
by K\"uhnau. Yet their exact values are known only
in a few special cases. Previously, the class of planar domains with
polygonal boundary has been studied by many authors from the
point of view of numerical computation. We present here a new
$hp$-FEM algorithm for the computation of moduli of rings
and quadrilaterals and compare its accuracy and performance
with previously known methods such as the Schwarz-Christoffel
Toolbox of Driscoll and Trefethen. We also demonstrate that the
$hp$-FEM algorithm applies to the case of non-polygonal
boundary and report results with concrete error bounds.
\end{abstract}



\pagestyle{myheadings}
\thispagestyle{plain}
\markboth{H. HAKULA, A. RASILA AND M. VUORINEN}{ON MODULI OF RINGS AND QUADRILATERALS}

\section{Introduction}

Plane domains with piecewise-smooth boundary curves occur in applications to electronics circuit design, airfoil modelling in computational fluid dynamics, computer vision and various other problems of engineering and science \cite{hen,kuh,kyt,ps,sl,sm}.
We assume that the domain is bounded and that there
are either one or two simple (and nonintersecting) boundary curves.
The domain is then either simply or doubly connected. For the mathematical
modelling of these domains it is usually convenient to map the domains
conformally onto ``canonical domains'' as simple as possible: the unit
disk $ \BC=\{z \in \C: |z|<1\}$ or the annulus
$ \{z \in \C: r<|z|<1\}\,.$ Sometimes a rectangle is
preferable to the unit disk as a canonical domain. The existence of
these canonical conformal mappings is guaranteed by classical results of
geometric function theory but the construction of this mapping in a concrete
application case is usually impossible. Therefore one has to resort
to numerical conformal mapping methods for which there exists an
extensive literature \cite{DrTr,kyt,ps2,sl}. The Schwarz-Christoffel (SC) Toolbox of
Driscoll \cite{Dr}, based on the software of Trefethen \cite{t}, is
in wide use for numerical conformal mapping applications.

In the doubly connected case, one might be interested in only knowing the inner
radius $r$ of the canonical annulus. For instance this occurs if we wish
to compute the electric resistance of a ring condenser. In this situation
the conformal mapping itself is not needed if we are able to find the
inner radius $r$ by some other method. It is a classical fact
that the inner radius $r$ can be obtained in terms of
the solution of the Dirichl\'et problem for the Laplace equation
in the original domain with the boundary value $0$ on one boundary component
and the boundary value $1$ on the other one. This fact is just one way of
saying that the modulus of a ring domain is conformally invariant: for the
canonical annulus $ \{z \in \C: r<|z|<1\}\,$ the modulus is
equal to $\log (1/r)\,.$ This idea reduces the problem of computing the
number $r$ to the problem of numerical approximation of the solutions
of Laplace equation in ring domains. In the paper \cite{bsv} this method
was applied to several concrete examples of ring domains for which
numerical results were reported, too. Again, it is also possible to use
the Schwarz-Christoffel method for doubly connected domains \cite{Hu}.

We next consider the case of simply connected plane domains. For such a
domain $D$ and for a quadruple $\{ z_1, z_2, z_3, z_4 \}$ of its boundary
points we call $(D;z_1, z_2,z_3,z_4)$ a quadrilateral if
$z_1, z_2,z_3,z_4$ occur in this order when the boundary curve is
traversed in the positive direction. The points $z_k,k=1,..,4\,,$
are called the vertices and the part of the oriented boundary between
two successive vertices such as $z_1$ and $z_2$ is called a boundary
arc $(z_1, z_2)\,.$ The modulus
$\M(D;z_1, z_2,z_3,z_4)$ of the quadrilateral $(D;z_1, z_2,z_3,z_4)$
is defined to be the unique $h>0$ for which there exists a conformal
mapping of $D$ onto the rectangle with vertices $1+ih,ih,0,1$ such that
the points $z_1,z_2,z_3,z_4$ correspond to the vertices in this order.
This conformal mapping is called the canonical conformal mapping associated
with the quadrilateral. As in the case of doubly connected domains
discussed above, it is well-known that the computation of the modulus $h$
of the quadrilateral may be reduced to solving the Dirichl\'et-Neumann
boundary value problem in the original domain $D$ with
the Dirichl\'et boundary values
$1$ on the boundary arc $(z_1,z_2)$ and $0$ on the arc $(z_3,z_4)$ and
Neumann boundary values $0$ on the arcs $(z_3,z_4)$ and $(z_4,z_1)\,.$

Conformal moduli of rings and quadrilaterals have independent theoretical
interest because of their crucial role in the theory of quasiconformal mappings \cite{LV}.
These quantities are closely related to certain physical constants, e.g.
they play an important role in applications involving the measurement of resistance
values of integrated circuit networks. But the problem of computing the moduli is also
interesting in the wider engineering context. The reciprocal identities (\ref{recipridty}) and (\ref{myidty})
can be used to generate test cases for curvilinear Lipschitz domains
and thus should be standard tools in the FEM-software development
community. Unfortunately these identities are missing from the introductory FEM textbooks. Although the experimental results in this paper show that the reciprocal identities provide error estimates similar to the true error (in cases where the exact analytic result is known) more investigations are needed to properly study their applicability in other contexts. Even though our interest lies in the high-order methods, these test cases are equally valid for any numerical PDE methods and mesh adaptation in particular.

One specific application area of the algorithms presented here is
the simulation of measurements for the 2D electrical impedance tomography (EIT) \cite{nh}.
In EIT problems a number of electrodes are placed on the boundary of the domain
and current patterns are considered between every pair of them. Indeed, computing
the moduli can be considered as a very crude model for the so-called  EIT background
forward problem. In general, the meshes for the EIT forward problems can be adapted using
the approaches outlined below. High level of accuracy is necessary for precise
control of artificial noise in the simulations.

A general observation about the literature seems to
be that reported numerical values of the moduli of concrete
quadrilaterals (or ring domains) are hard to find.
Perhaps the longest list of numerical results is
given in \cite{bsv} where pointers to earlier literature may be found.
The recent book \cite{ps2} lists also many such numerical values.
In our opinion a catalogue of these numerical values in the simplest
cases would be desirable for instance for reference purposes.
The book \cite{sl} and the paper \cite[p. 127]{ps} list certain
engineering formulas which have been applied in VLSI circuit design.

An outline of the structure of this paper now follows. First, in Section 2
we describe the methods used in this paper. In Section 3 we discuss in detail
the various FEM methods used here, in particular the $hp$-method which was
implemented and applied to generate some of the results reported below.
Another method we use is the $h$-adaptive software package
AFEM of K. Samuelsson, which implements an adaptive FEM method and which
was previously used in \cite{bsv}. In the present paper we use the AFEM
method to compute the modulus of a quadrilateral whereas in \cite{bsv}
it was used merely for the computation of the moduli of ring domains.
In Section 4 a test problem for quadrilaterals is described together
with its analytic solution, following \cite{hvv}. This analytic solution
requires, however, an application of a numerical root finding program.
Accordingly, this formula is analytic-numeric in its character.
In Section 5 we check several methods against this analytic formula in
a test involving a family of convex quadrilaterals.
The methods discussed are the analytic
formula from \cite{hvv}, the Schwarz-Christoffel Toolbox of \cite{Dr,DrTr},
the AFEM method of Samuelsson \cite{bsv} and the present $hp$-method.
On the basis of these experiments, an accuracy ranking
of the methods is given in Section 5.
In Section 6 the more general case of polygonal quadrilaterals is
investigated, in particular L-shaped domains, and the results are compared
to the literature. In Section 7 we discuss the computation of the  modulus of a ring domain in a few special cases. For instance, for ``the cross in square'' ring domain considered previously in \cite[Example 4]{bsv} we now obtain much improved accuracy. In Section 8 we compute some examples with the $hp$-FEM which are difficult for other methods. In Section 9 our results and discoveries are summarized.

\section{Methods}

The following problem is known as the {\it Dirichl\'et-Neumann problem}.
Let $D$ be a region in the complex plane whose boundary
$\partial D$ consists of a finite number of regular Jordan
curves, so that at every point, except possibly at finitely many points,
of the boundary a normal is defined. Let $\partial D =A \cup B$
where $A, B$ both are unions of Jordan arcs. Let $\psi_A,\psi_B$ be a real-valued
continuous functions defined on $A,B$, respectively. Find a function $u$
satisfying the following conditions:
\begin{romannum}
\item
$u$ is continuous and differentiable in
$\overline{D}$.
\item
$u(t) = \psi_A(t),\qquad \textrm{ for all}\,\, t \in A$.
\item
If $\partial/\partial n$ denotes differentiation in
the direction of the exterior normal, then
$$
\frac{\partial}{\partial n} u(t)=\psi_B(t),\qquad \textrm{ for all}\,\, t\in B.
$$
\end{romannum}

\subsection{Modulus of a quadrilateral and Dirichl\'et integrals}
One can express the modulus of a quadrilateral $(D; z_1, z_2, z_3, z_4)$
in terms of the solution of the Dirichl\'et-Neumann problem as follows.
Let $\gamma_j, j=1,2,3,4$ be the arcs of
$\partial D$ between $(z_4, z_1)\,,$ $(z_1, z_2)\,,$ $(z_2, z_3)\,,$
$(z_3, z_4),$ respectively. If $u$ is the (unique) harmonic solution of
the Dirichl\'et-Neumann problem with boundary values of $u$ equal to $0$ on
$\gamma_2$, equal to $1$ on $\gamma_4$ and with $\partial u/\partial n =
0$ on $\gamma_1 \cup \gamma_3\,,$ then by \cite[p. 65/Thm 4.5]{Ah}:
\begin{equation}
\label{qmod}
\M(D;z_1,z_2,z_3,z_4)=
\int_D |\nabla
u|^2\,dm.
\end{equation}

\subsection{Modulus of a ring domain and Dirichl\'et integrals}
Let $E$ and $F$ be two disjoint compact sets in the extended
complex plane ${\C_\infty}$. Then one of the sets $E,$ $F$ is bounded and
without loss of generality we may assume that it is $E\,.$ If both $E$ and $F$ are connected
and the set $R={\C_\infty} \setminus(E \cup F)$ is connected,
then $R$ is called a {\it ring domain}. In this case $R$ is a doubly
connected plane domain. The {\it capacity} of $R$ is defined by
$$
\Cap R=\inf_u \int_R |\nabla u|^2\,dm,
$$
where the infimum is taken over all nonnegative, piecewise
differentiable functions $u$
with compact support in $R\cup E$ such that $u=1$ on $E$.
It is well-known that the harmonic function on $R$ with
boundary values $1$ on $E$ and $0$ on $F$ is the unique function that
minimizes the above integral. In other words, the minimizer may be found
by solving the Dirichl\'et problem for the Laplace equation in $R$ with
boundary values $1$ on the bounded boundary component $E$ and $0$ on the
other boundary component $ F\,.$
A ring domain $R$ can be
mapped conformally onto the annulus $\{z:e^{-M}<|z|<1\}$, where
$M=\M(R)$ is the {\it conformal modulus}
of the ring domain $R\,.$ The modulus and capacity of a ring
domain are connected by the simple identity
$\M(R)=2\pi/\Cap R$.
For more information on the modulus of a ring domain
and its applications in complex analysis the reader is
referred to \cite{Ah,hen,kuh,LV,ps2}.

\subsection{Classification of methods for numerical computing}
For the computation of the modulus of a quadrilateral or
of a ring domain there are two natural approaches
\begin{romannum}
\item methods based on the definition of the modulus and on the
use of a conformal mapping onto a canonical rectangle or annulus,
\item methods that give only the modulus, not the canonical
conformal map.
\end{romannum}
In some sense, methods of class (i) give a lot of extra information,
namely the conformal mapping -- all we want is a single real number.
Methods of class (ii) rely on solving the Dirichl\'et-Neumann boundary
value problem or Dirichl\'et problem for the Laplace equation
as described above.

In this paper we will mainly use methods of type (ii) that make use of
adaptive FEM methods for solving the Laplace equation.

\subsection{Review of the literature on numerical conformal mapping}
With the exception of a few special cases, both of the above methods lead to
extensive numerical computation. For both classes of methods
there are several options in the literature, see for instance the
bibliography of \cite{bsv}. Various aspects of the theory and practice
of numerical conformal mapping are reviewed in the monographs
\cite{DrTr,kyt,ps2,sl}.
See also the authoritative surveys \cite{gai,Pap,td,weg}.

Recently numerical conformal mappings have been studied from
various points of view and in various applications by many authors,
see e.g. \cite{a,b,c,cm,ddep,k,mr,p1,Porter}. In \cite[Chapter 3]{ps2} N.~Papamichael and N.~Stylianopoulos describe the so-called domain decomposition method for the computation of the modulus of a quadrilateral which is designed for the case of elongated quadrilaterals and applies e.g. to polygonal quadrilaterals that can be decomposed into simple pieces whose moduli
can be estimated. As an example they consider a spiraling  quadrilateral that can be decomposed into a ``sum'' of $13$ trapezoids and report results that are expected to be correct up to $7$ decimal places. Therefore, this method seems very attractive for the computation of the modulus of a special class of quadrilaterals. A key feature of the method is that it reduces the numerical difficulties caused by the crowding phenomenon for this special class of quadrilaterals.

\section{$p$-, and $hp$-finite element method}

In the paper \cite{bsv} the modulus of a ring domain was computed with the help of the software package AFEM of K.~Samuelsson, based on an $h$-adaptive finite element method. It can be easily applied to compute the modulus of a quadrilateral.

In this section we describe the high-order $p$-, and $hp$-finite element methods. The paper of Babu\v{s}ka and Suri \cite{bs} gives an accessible overview of the method. For a more detailed exposition we refer to Schwab \cite{s}, and for those familiar with engineering approach the book by Szabo and Babu\v{s}ka \cite{sz} is recommended. For the applications considered in this paper, any finite element computation requires at least the choice of the following.
\begin{remunerate}
\item Initial discretization of the domain. In 2D each discretization or mesh divides the domain into elements, plane regions with piecewise smooth boundaries. These are usually either triangles or quadrilaterals.
\item Refinement strategy. The choice of the refinement strategy is connected
to choosing the finite element method (FEM): mesh refinement ($h$-method), elementwise polynomial order ($p$-method), or both above ($hp$-method). The unknowns or degrees of freedom are the coefficients of the chosen shape functions. In the $h$-version the shape functions are such that the coefficients are also values of the solution at specified locations of the discretization of the computational domain, that is, the nodes of the mesh. In the $p$-method, the shape functions are polynomials that are associated with topological entities of the elements, nodes, sides, and interior. Thus, in addition to increasing accuracy through refining the mesh, we have an
additional refinement parameter, the polynomial degree $p$.
\end{remunerate}

Both choices will have an influence on the performance and the accuracy attained with the chosen method. The mutual influence of these choices is hard to analyze theoretically but usually it may be seen in the results. For instance, we have observed that the choice of the intitial mesh and the mesh refinement strategy may limit the accuracy achieved by the $hp$-method and therefore it is useful to try a few initial meshes.

Let us next define a $p$-type quadrilateral element. The construction of triangles
is similar and can be found from the references given above.

\subsection{Shape functions}
Many different selections of shape functions are possible.
We follow Szabo and Babu\v{s}ka \cite{sz} and present the so-called
hierarchic shape functions.

Legendre polynomials of degree $n$ can be defined using a recursion formula
\begin{equation}
(n+1)P_{n+1}(x) - (2n+1) x P_n(x) + n P_{n-1}(x) =0,\quad P_0(x)=1.
\end{equation}
The derivatives can similarly be computed using a recursion
\begin{equation}
(1-x^2)P'_n(x)=-n x P_n(x) + n P_{n-1}(x).
\end{equation}

For our purposes the central polynomials are the integrated Legendre polynomials for $x\in[-1,1]$
\begin{equation}
\phi_n(\xi)=\sqrt{\frac{2n-1}{2}}\int_{-1}^{\xi}P_{n-1}(t)\,dt,\quad n=2,3,\ldots
\end{equation}
which can be rewritten as linear combinations of Legendre polynomials
\begin{equation}
\phi_n(\xi)=
\frac{1}{\sqrt{2(2n -1)}}\left(P_n(\xi) - P_{n-2}(\xi)\right),\quad n=2,3,\ldots
\end{equation}
The normalizing coefficients are chosen so that
\begin{equation}
\int_{-1}^{1}\frac{d \phi_i(\xi)}{d\xi} \frac{d \phi_j(\xi)}{d\xi} \,d\xi =
	\delta_{ij}, \quad i,j \geq 2.
\end{equation}

We can now define the shape functions for a quadrilateral reference element
over the domain $[-1,1]\times[-1,1]$. The shape functions are divided into three categories:
nodal shape functions, side modes, and internal modes.

\subsection{Nodal shape functions} There are four nodal shape functions:
\begin{align*}
N_1(\xi,\eta) &= \frac{1}{4}(1-\xi)(1-\eta), \\
N_2(\xi,\eta) &= \frac{1}{4}(1+\xi)(1-\eta), \\
N_3(\xi,\eta) &= \frac{1}{4}(1+\xi)(1+\eta), \\
N_4(\xi,\eta) &= \frac{1}{4}(1-\xi)(1+\eta). \\
\end{align*}
Taken alone, these shapes define the standard four-node quadrilateral
finite element.

\subsection{Side shape functions} There are $4(p-1)$ side modes associated with
the sides of a quadrilateral $(p\geq 2)$.
\begin{align*}
N_i^{(1)}(\xi,\eta) &= \frac{1}{2} (1-\eta) \phi_i(\xi), \quad i=2,\ldots,p,\\
N_i^{(2)}(\xi,\eta) &= \frac{1}{2} (1+\xi) \phi_i(\eta),\quad i=2,\ldots,p,\\
N_i^{(3)}(\xi,\eta) &= \frac{1}{2} (1+\eta) \phi_i(\eta),\quad i=2,\ldots,p,\\
N_i^{(4)}(\xi,\eta) &= \frac{1}{2} (1-\xi) \phi_i(\xi), \quad i=2,\ldots,p.\\
\end{align*}

\subsection{Internal shape functions} For the internal modes we have two options.
The so-called trunk space has $(p-2)(p-3)/2$ shapes
\begin{equation}
N_{i,j}^0(\xi,\eta)=\phi_i(\xi)\phi_j(\eta),\quad i,j\geq 2,\quad i+j=4,5,\ldots,p,
\end{equation}
whereas the full space has $(p-1)(p-1)$ shapes
\begin{equation}
N_{i,j}^0(\xi,\eta)=\phi_i(\xi)\phi_j(\eta),\quad i=2,\ldots,p,\quad j=2,\ldots,p.
\end{equation}
In this paper we always use the full space.
The internal shape functions are often referred to as bubble-functions.

\subsection{Parity problem}
The Legendre polynomials have the property $P_n(-x)=(-1)^n P_n(x)$.
In 2D all internal edges of the mesh are shared by two different elements.
We must ensure that each edge has the same global parameterization in
both elements. This additional book-keeping is not necessary in the standard
$h$-FEM.

\subsection{Resource requirements}
We have seen that the number of unknowns in a $p$-type quadrilateral is
$(p+2)(p+3)/2-p$ or $4p+(p-1)^2$ if the internal modes are from
trunk or full space, respectively.
To compensate this, the number of elements is naturally taken to be as small as
possible.
Indeed, when the mesh is adapted in a suitable way, the dimension of the
overall linear system can be significantly lower than in the corresponding $h$-method.
However, the matrices tend to be denser in the $p$-method, so
the space requirements in relation to the dimension of the linear system are
greater for the $p$-method.

\subsection{Proper grading of the meshes}
For a certain class of problems it can be shown that if the mesh and the elemental
degrees have been set optimally, we can obtain \textit{exponential convergence}.
A geometric mesh is such that each successive layer of elements changes in
size with some \textit{geometric scaling factor} $\alpha$,
toward some point of interest. In this case, the points of interest are always corner points.

The following theorem is due to Babu{\v{s}}ka and Guo \cite{bg1}. Note that construction of
appropriate spaces is technical. For rigorous treatment of the theory involved
see Schwab \cite{s},
Babu{\v{s}}ka and Guo \cite{bg2} and references therein.

\begin{theorem} \label{propermesh}
Let $\Omega \subset \mathbb{R}^2$ be a polygon, $v$ the FEM-solution, and
let the weak solution $u_0$ be in a suitable countably normed space where
the derivatives of arbitrarily high order are controlled.
Then
\[
\inf_v \|u_0 - v\|_{H^1(\Omega)} \leq C\,\exp(-b \sqrt[3]{N}),
\]
where $C$ and $b$ are independent of $N$, the number of degrees of freedom.
Here $v$ is computed on a proper geometric mesh, where the orders of individual
elements depend on their originating layer, such that the highest layers have the smallest
orders.

The result also holds for constant polynomial degree distribution.
\end{theorem}

Let us denote the number of the highest layer with $\nu$, \textit{the nesting level}.
Using this notation we can refer to geometric meshes as $(\alpha,\nu)$-meshes.

\begin{figure}
\begin{center}
\includegraphics[width=.35\textwidth]{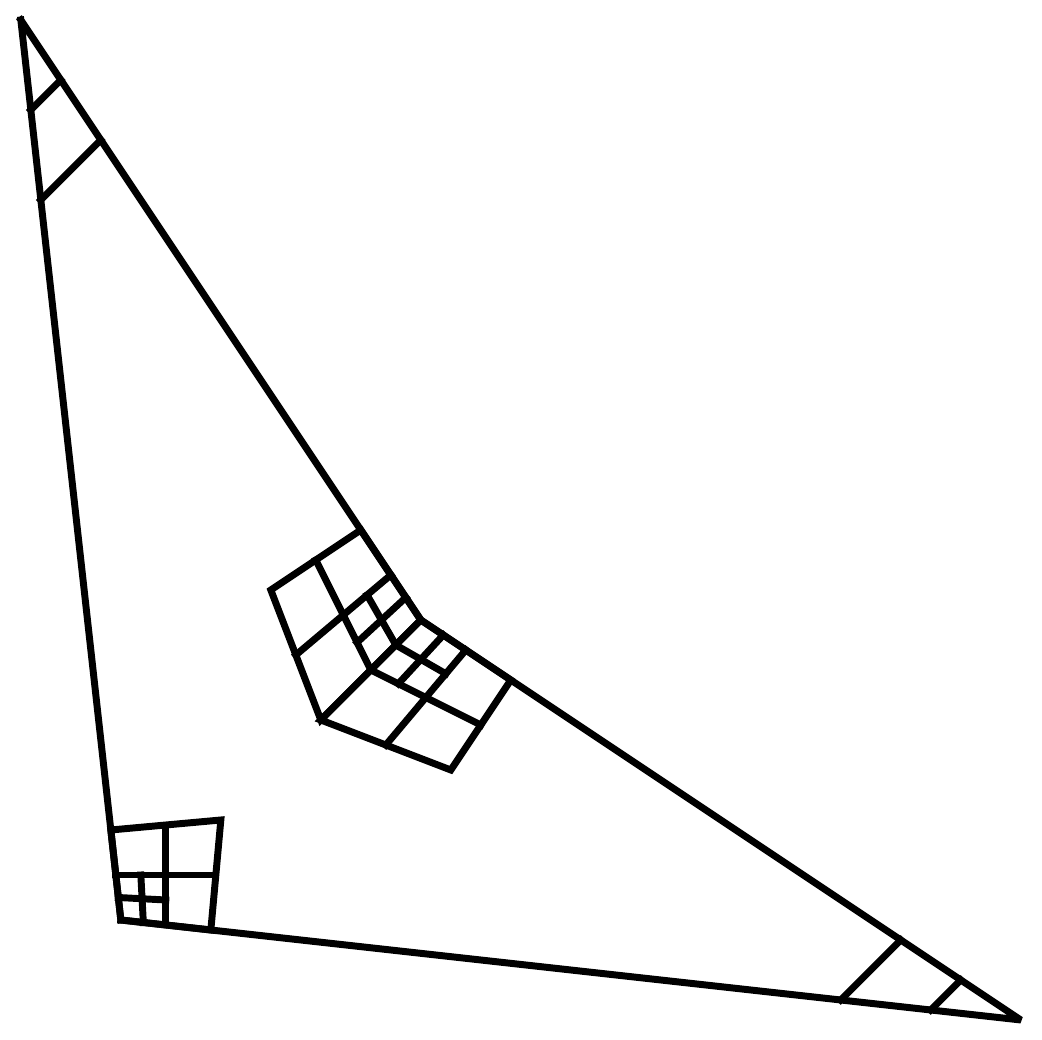}
\end{center}
\caption{Geometric mesh for a general quadrilateral.}\label{geometricmeshfig}
\end{figure}

In Figure \ref{geometricmeshfig} we show a geometric mesh template for
a non-convex quadrilateral. Here we require that each node lies at the end
point of an edge and that the meshlines follow the guidelines of
the geometric meshes.

\begin{figure}
\begin{center}
\includegraphics[width=.25\textwidth]{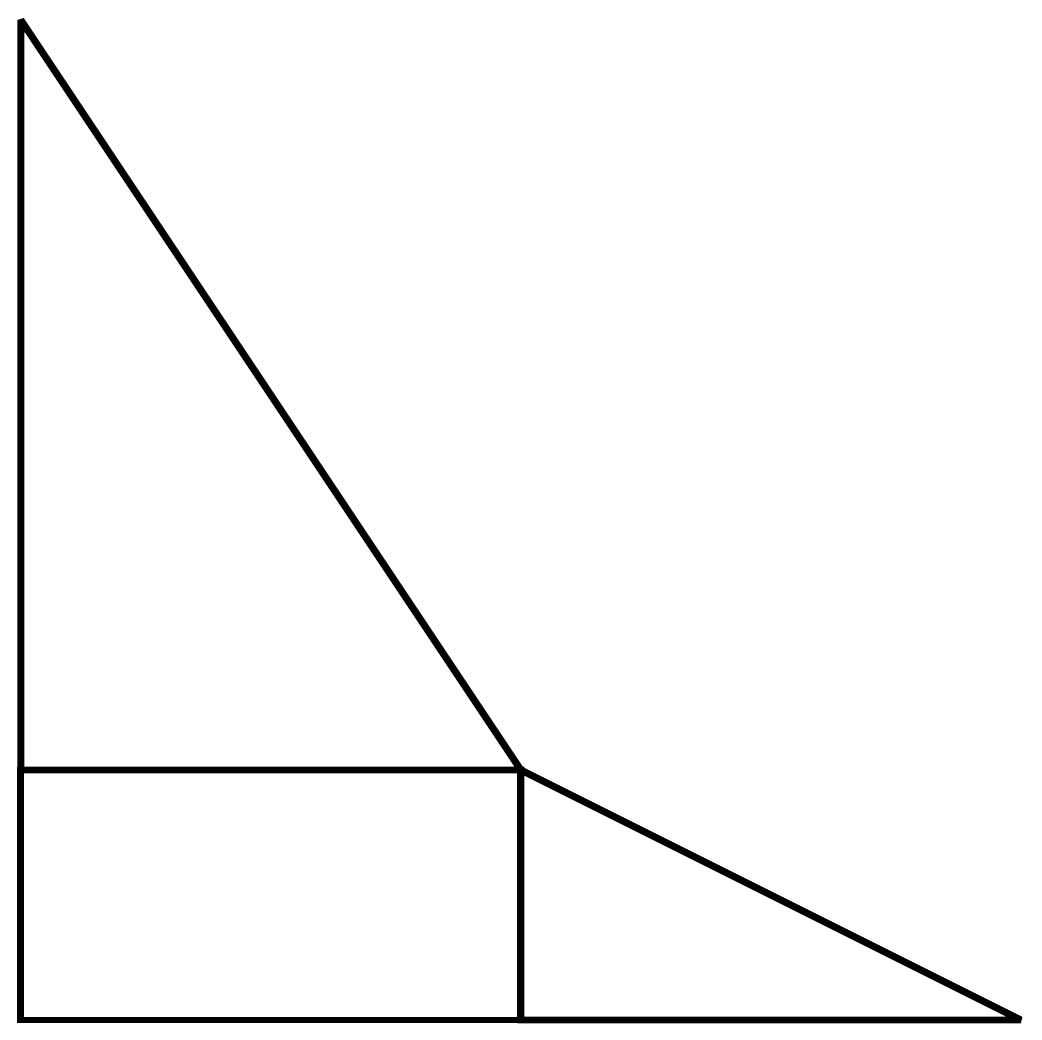}
\includegraphics[width=.25\textwidth]{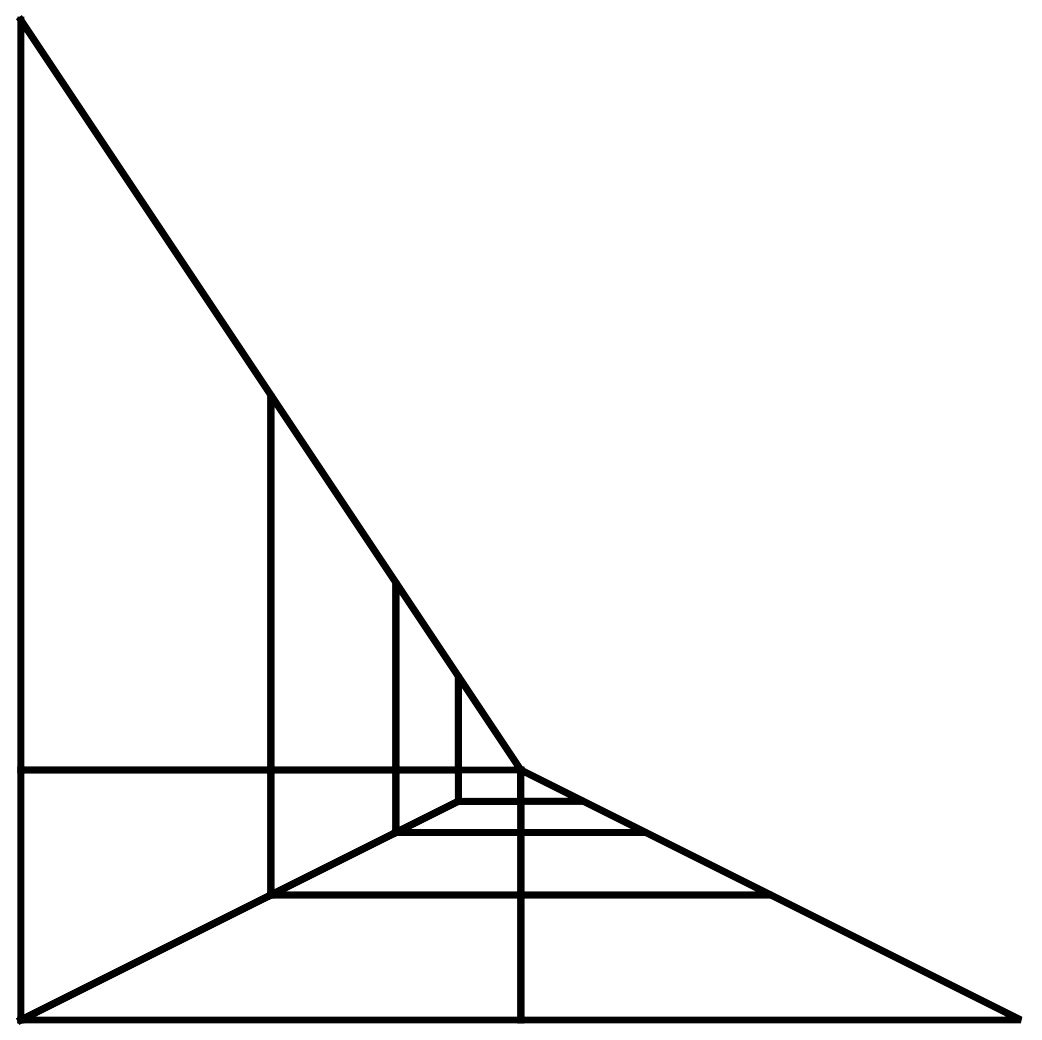}
\includegraphics[width=.25\textwidth]{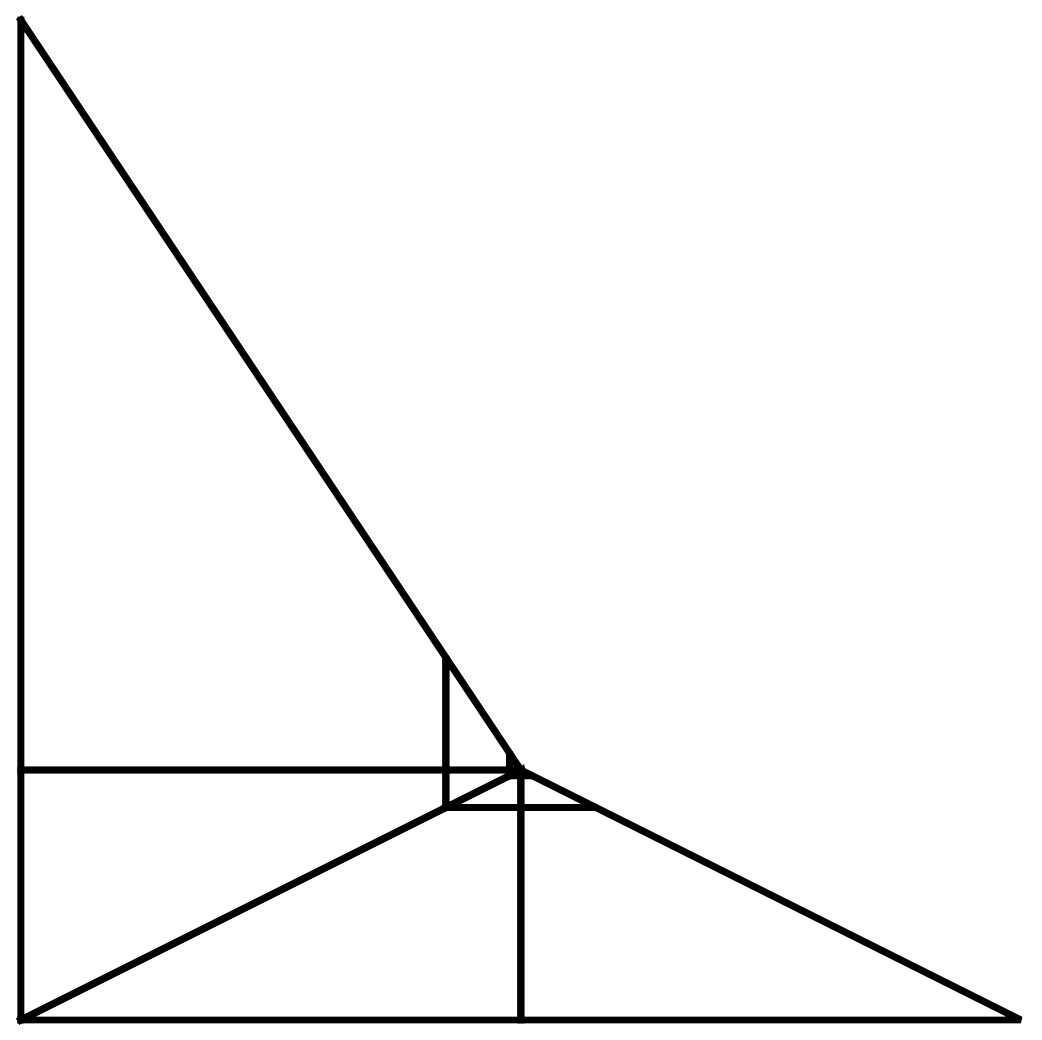}
\end{center}
\caption{Graded meshes: Effect of the scaling factor. From left to right,
template mesh, $(\alpha,\nu) = (1/2,3)$, $(\alpha,\nu) = (1/6,3)$.}\label{propermeshseqfig}
\end{figure}

In Figure \ref{propermeshseqfig} a sequence of graded meshes is shown.
In the middle and the rightmost meshes the number of elements is the same
because the nesting level is the same, only the scaling factor changes.

\subsection{Generating geometric meshes} \label{meshing}
Here we consider generation of geometric meshes in polygonal domains.
We use the following two-phase algorithm:
\begin{remunerate}
\item Generate a minimal mesh (triangulation) where the corners are
isolated with a fixed number of triangles depending on the interior angle, $\theta$ so that the refinements can be carried out independently:
\begin{romannum}
\item $\theta \leq \pi/2$: one triangle,
\item $\pi/2 < \theta \leq \pi$: two triangles, and
\item $\pi < \theta$: three triangles.
\end{romannum}
\item Every triangle attached to a corner is replaced by a refinement,
where the edges incident to the corner are split as specified by the scaling factor
$\alpha$. This process is repeated recursively until the desired nesting level $\nu$
is reached. Note that the mesh may include quadrilaterals after refinement.
\end{remunerate}

\begin{figure}
\begin{center}
\includegraphics[width=.25\textwidth]{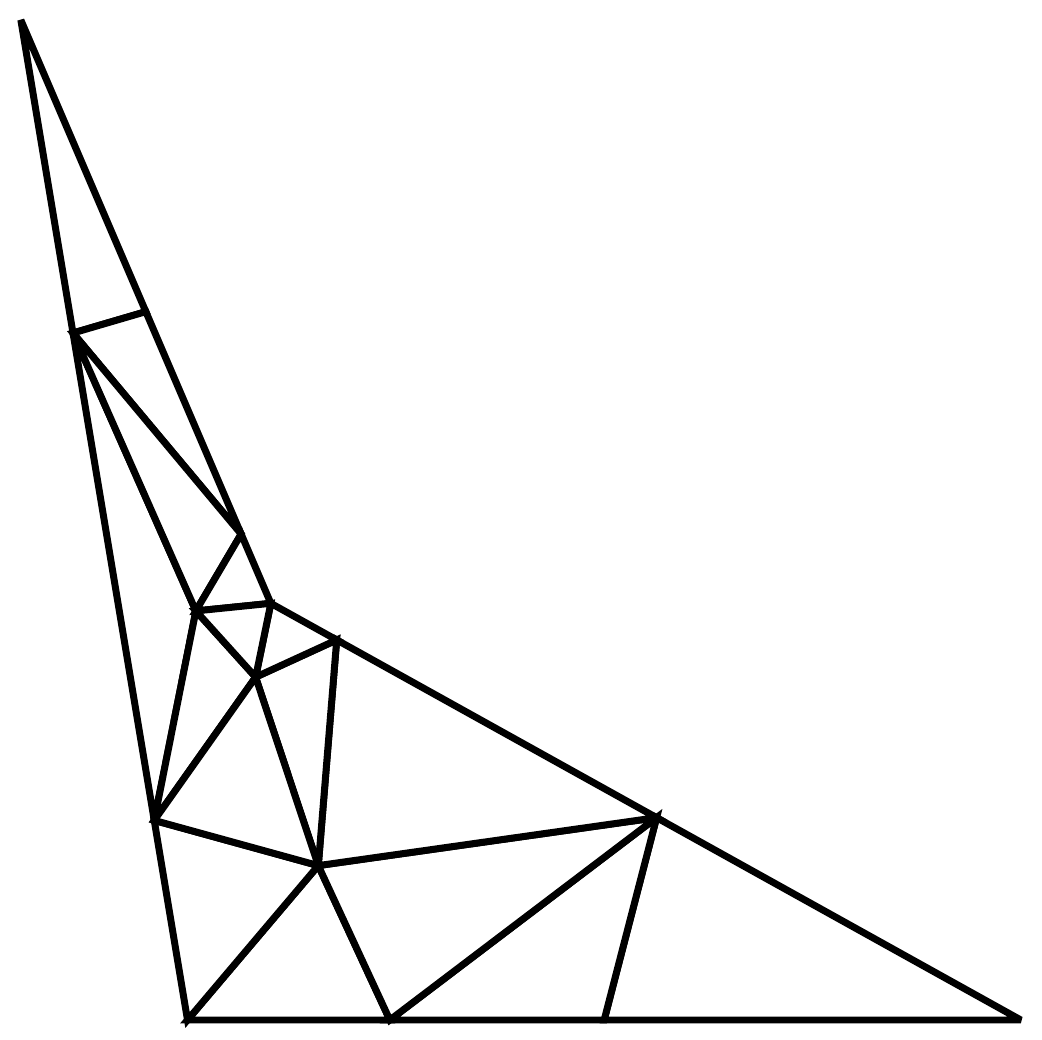}
\includegraphics[width=.25\textwidth]{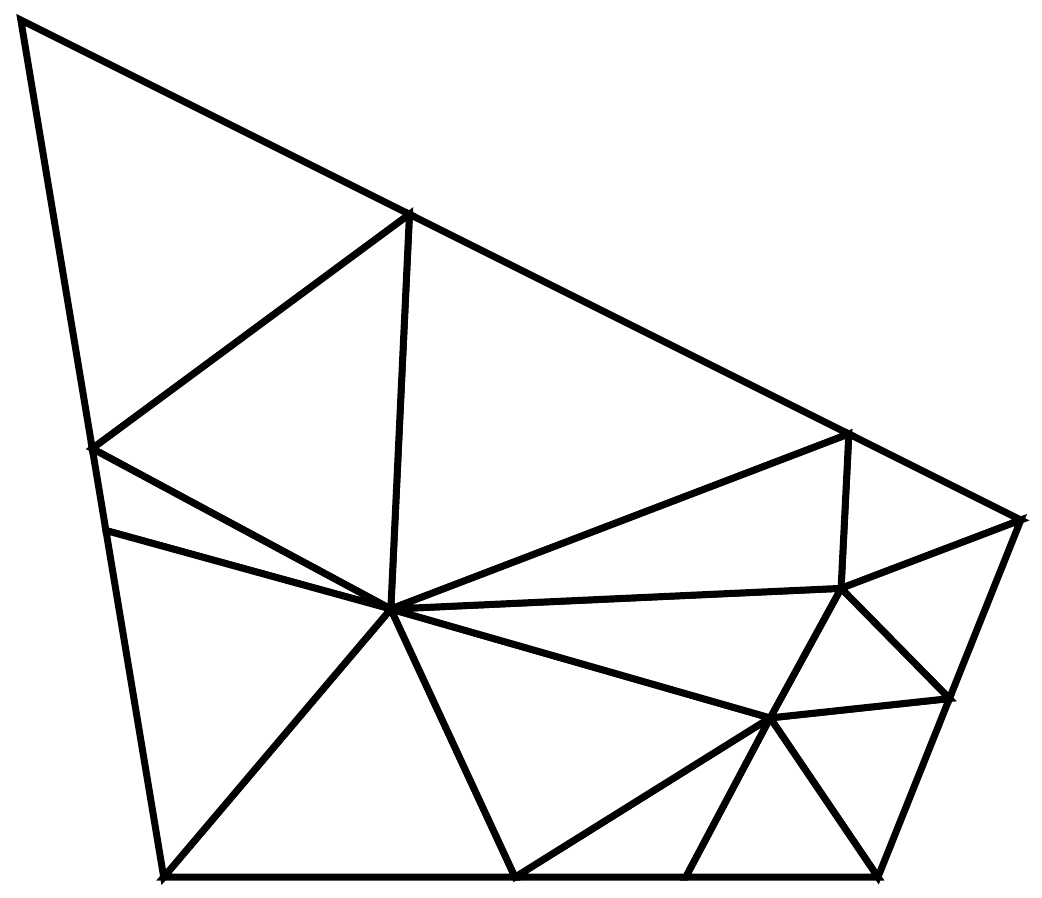}
\includegraphics[width=.25\textwidth]{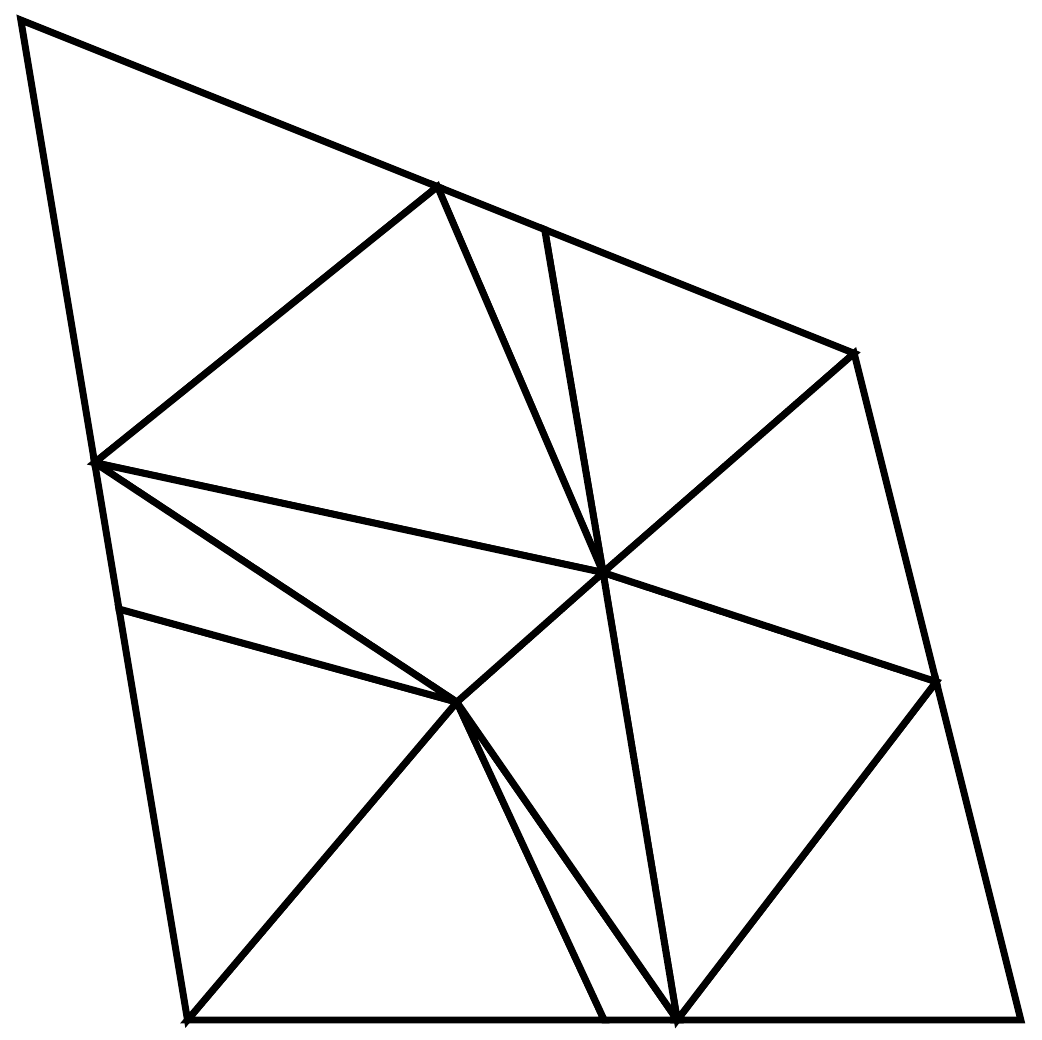}
\end{center}
\caption{Three sample meshes used in numerical experiments below. Note the refinement of the mesh structure close to the corner points.}\label{generatedmeshesfig}
\end{figure}

\begin{figure}
\begin{center}
\includegraphics[width=.31\textwidth]{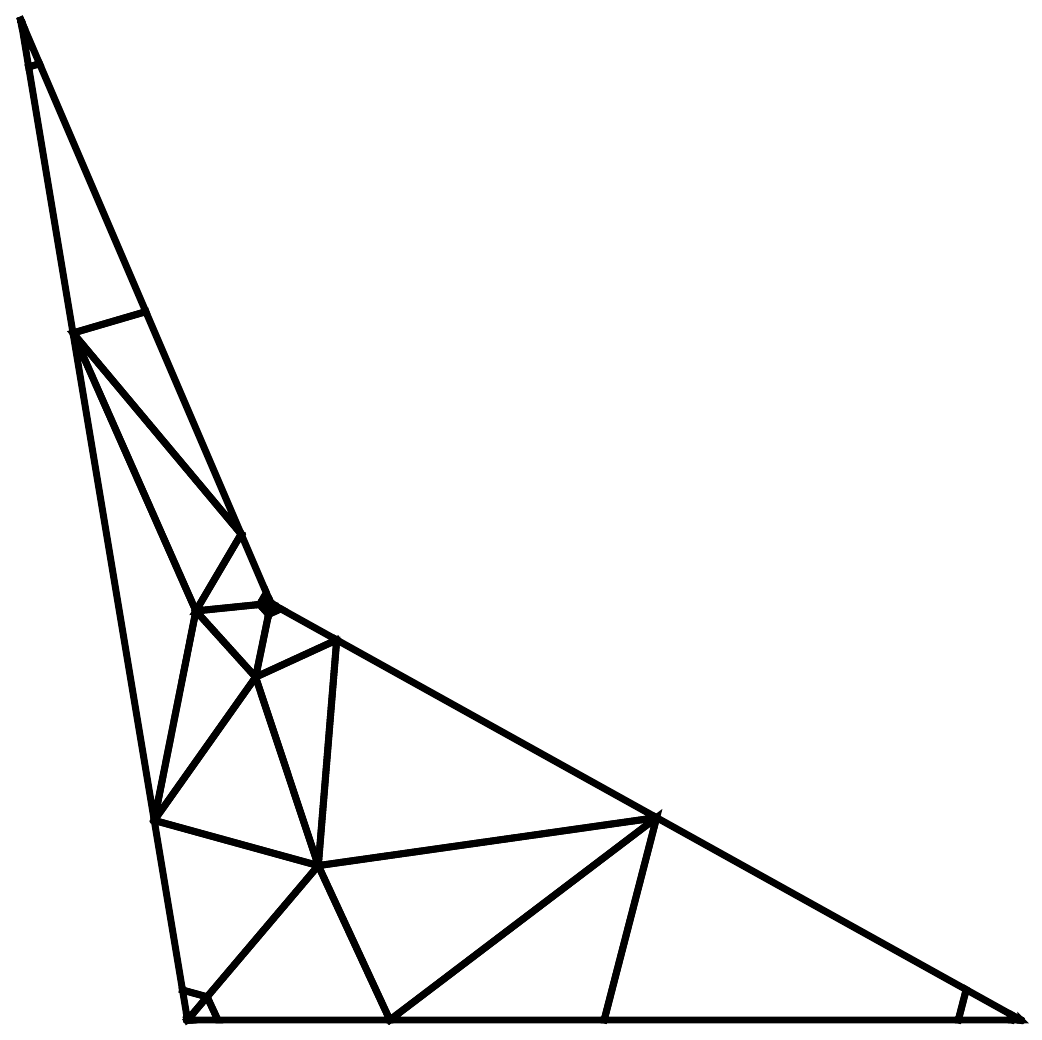}
\end{center}
\caption{Final geometric or $(0.15,12)$-mesh.
Due to small $\alpha$ only first two levels are visible.}
\label{generatedgradedmeshesfig}
\end{figure}

In Figure \ref{propermeshseqfig} we can also see our preferred
element subdivisions: triangle to (quadrilateral, triangle)-pair, and
quadrilateral to three quadrilaterals. These two rules are sufficient for our purposes since
we always grade toward a corner point. Using this, we can derive a simple estimate for the number of degrees of freedom
$N$. Letting $T$ denote the number of elements in the initial mesh, and $C$ the number of corners in the domain
(or those used in refining):
\begin{equation}\label{dofestimate}
N \sim (T + 6C \nu)p^2,
\end{equation}
where the constant 6 is the product of the maximal number of elements surrounding a corner, 3,
and the maximal number of new elements per level, 2.

Finally, in Figure \ref{generatedmeshesfig} three minimal meshes and
in Figure \ref{generatedgradedmeshesfig} one final mesh are shown.

\begin{figure}
\begin{center}
\includegraphics[width=0.6\textwidth]{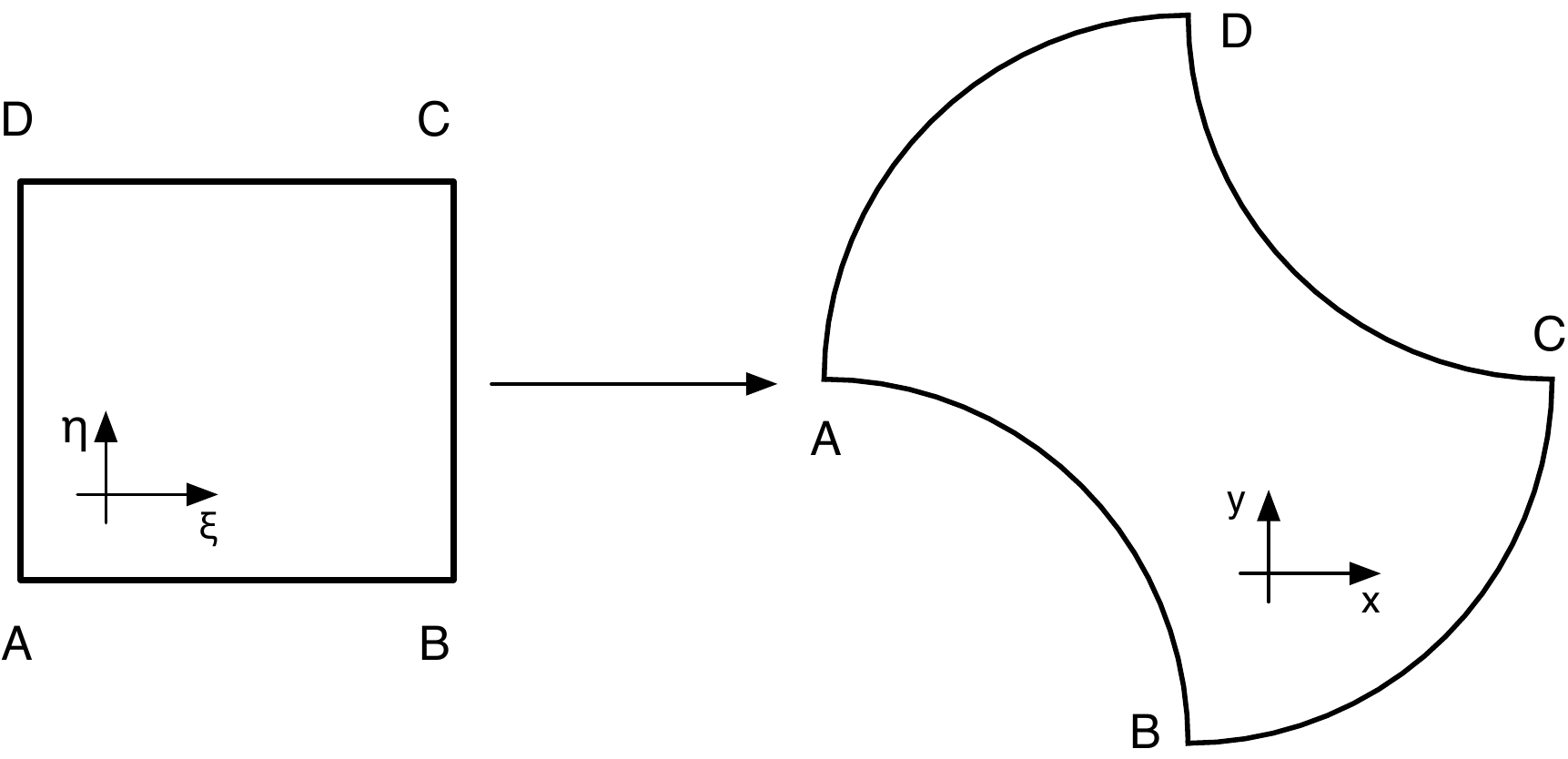}
\end{center}
\caption{Curved boundary mapping.}\label{curvedboundary}
\end{figure}

\subsection{Domains with curved boundaries}
Since we want to use as large elements as possible, it is important to represent
curved boundary segments accurately. The linear blending function method
of Gordon and Hall \cite{gh} is our choice for this purpose.

In the general case all sides of an element can be curved as in Figure~\ref{curvedboundary}.
We assume that every side is parameterized:
\begin{equation}
x=x_i(t),\ y=y_i(t),\quad -1\leq t \leq 1,\quad i=1,\ldots,\mbox{number of sides}
\end{equation}
Using capital letters as coordinates of the corner points, $(X_i,Y_i)$, we can write the mapping
for the global $x$-coordinates of a quadrilateral as
\begin{equation}
\begin{split}
x& = \frac12(1-\eta)x_1(\xi)+\frac12(1+\xi)x_2(\eta)+\frac12(1+\eta)x_3(\xi)+\frac12(1-\xi)x_4(\eta)\\
&-\frac14(1-\xi)(1-\eta)X_1-\frac14(1+\xi)(1+\eta)X_2-\frac14(1+\xi)(1+\eta)X_3\\
&-\frac14(1-\xi)(1+\eta)X_4,
\end{split}
\end{equation}
and symmetrically for the $y$-coordinate. Note, that if the side parameterizations represent
straight edges, the mapping simplifies to the standard bilinear mapping of quadrilaterals.

In the following we always use exact representation of the geometry which implies that
in the ensuing mesh grading process no approximation of geometry is necessary.
Here the mesh generation of the curved domains is template-based, thus the changes
in curvature are not automatically dealt with. For a highly accessible review of the $p$-method
mesh generation issues we refer to \cite{pm}.

\section{Convex quadrilateral}

In this section our goal is to introduce a test problem, whose
solution is determined by a transcendental equation.
This equation can be numerically solved to the desired
accuracy and we will use this to check the
validity of the numerical methods we use as well as to obtain an
experimental estimate for their accuracy. The test problems we consider
are convex polygonal quadrilaterals.
The simplest such quadrilateral consists of the four vertices
and the line segments joining the vertices. Let $z_1,z_2,z_3,z_4
\in \mathbb{C}$ be distinct points and suppose that the polygonal
line that results from connecting these points by segments in the
order $z_1,z_2,z_3,z_4, z_1$ forms the positively oriented boundary
of a domain $Q$. For simplicity, we denote by $\QM(z_1,z_2,z_3,z_4)$
the modulus $\M(Q; z_1,z_2,z_3,z_4)$. Then the modulus is a conformal
invariant in the following sense: If $f\colon Q \to fQ$ is a conformal
mapping onto a Jordan domain $fQ,$ then $f$ has a homeomorphic extension
to the closure $\overline{Q}$ (also denoted by $f$) and
$$
\M(Q; z_1,z_2,z_3,z_4) = \M(fQ; f(z_1),f(z_2),f(z_3),f(z_4)) \,.
$$
It is clear by geometry that the following reciprocal identity holds:
\begin{equation} \label{recipridty}
\M(Q; z_1,z_2,z_3,z_4) \M(Q; z_2,z_3,z_4, z_1) =1.
\end{equation}

There are two particular cases, where we can immediately give
$\QM(z_1,z_2,z_3,z_4).$ The first case occurs, when all the sides
are of equal length (i.e. the quadrilateral is a rhombus) and in this
case the modulus is $1\,,$ see \cite{hvv}. In the second case
$(Q; z_1,z_2,z_3,z_4)$ is the rectangle $(Q; 1+ih,ih,0,1)$, $h>0$, and
$\QM(1+ih,ih,0,1)=h.$

\subsection{Basic identity}
In \cite[2.11]{hvv} some identities satisfied by the function
$\QM(a,b,0,1)$ were pointed out. We will need here the following
one, which is the basic reciprocal identity (\ref{recipridty})
rewritten for the expression $\QM$ :
\begin{equation} \label{myrecidty}
\QM(a,b,0,1)\cdot\QM((b-1)/(a-1),1/(1-a), 0,1)=1 \,.
\end{equation}

We shall consider here the following particular cases of this
reciprocal identity: (a) parallelogram, (b)
trapezoid with angles $(\pi/4,3\pi/4, \pi/2,\pi/2)$, and (c)
a convex polygonal quadrilateral.
Note that for the cases (a) and (b) the formula is less complex
than for the general case (c).

\subsection{The hypergeometric function and complete elliptic integrals}
Given complex numbers
$a,b,$ and $c$ with $c\neq 0,-1,-2, \ldots $,
the {\em Gaussian hypergeometric function} is the analytic
continuation to the slit plane $\C \setminus [1,\infty)$ of
the series
\begin{equation} \label{eq:hypdef}
F(a,b;c;z) = {}_2 F_1(a,b;c;z) =
\sum_{n=0}^{\infty} \frac{(a,n)(b,n)}{(c,n)} \frac{z^n}{n!}\,,\:\: |z|<1 \,.
\end{equation}
Here $(a,0)=1$ for $a \neq 0$, and $(a,n)$
is the {\em shifted factorial function}
or the {\em Appell symbol}
$$
(a,n) = a(a+1)(a+2) \cdots (a+n-1)
$$
for $n \in \N \setminus \{0\}$, where
$\N = \{ 0,1,2,\ldots\}$
and the
{\it elliptic integrals} $\K(r),\K'(r)$
are defined by
$$
\K(r)=\frac{\pi}{2} F(1/2,1/2;1; r^2),
\qquad
\K'(r)=\K(r'),\text{ and }r'= \sqrt{1-r^2}.
$$
Some basic properties of these functions can be found in \cite{avv} and \cite{OLBC}.

\subsection{Parallelogram}
For $t \in (0,\pi)$ and $h>0$ let
$$
g(t,h) \equiv \QM(1+ h e^{it}, h e^{it}, 0,1).
$$
An analytic expression for this function has been given in
\cite[2.3]{aqvv}:
\begin{equation} \label{aqvvform1}
g(t,h) = \K'(r_{t/\pi})/\K(r_{t/\pi}),
\end{equation}
where
\begin{equation} \label{aqvvform2}
r_a = \mu_a^{-1}\bigg(\frac{\pi h}{2\sin(\pi a)}\bigg),
\text{ for }0<a< 1,
\end{equation}
and the decreasing homeomorphism $\mu_a\colon (0,1)\to(0,\infty)$
is defined by
\begin{equation}
\label{mu_a}
\mu_a(r)\equiv \frac{\pi}{2\sin(\pi
a)}\,\frac{F(a,1-a;1;1-r^2)}{F(a,1-a; 1; r^2)}.
\end{equation}

\begin{theorem} \label{modquad} \cite{hvv}
Let $0 < a, b < 1$, $\max \{a+b,1 \} \le c \le 1 + \min \{a,b\}$,
and let $Q$ be the quadrilateral in the upper half plane $\HC=\{ z \in \C : \im z >0\}$
with vertices $0, 1, A$ and $B$, the interior angles at which are, respectively, $b\pi, (c-b)\pi,
(1-a)\pi$ and $(1+a-c)\pi$. Then the conformal modulus of $Q$
is given by
\begin{equation}
\label{hvvmodu}
\QM(A,B,0,1)\equiv \M(Q) = \K(r') / \K(r),
\end{equation}
where $r \in (0,1)$ satisfies the equation
\begin{equation} \label{eq:star}
A - 1 = \frac{L {r'}^{2(c-a-b)} F(c-a,c-b;c+1-a-b; {r'}^2)}{F(a,b;c;r^2)} \, ,
\end{equation}
say, and
$$
L = \frac{B(c-b,1-a)}{B(b,c-b)} e^{(b+1-c)i \pi}.
$$
\end{theorem}

For a fixed complex number $b$ with $\im(b)>0$ define the following function
$g(x,y)= \QM(x+i\cdot y,b,0,1)$ for $x\in \mathbb{R}$, $y>0\,.$ This is well-defined only if the polygonal domain with vertices $x+i\cdot y$, $b$, $0$, $1$ is positively oriented. This holds e.g. if $\re(b)<0$ and $x>0$.
It is a natural
question to study the level sets of the function $g\,.$ This function
tells us how the modulus of a polygonal quadrilateral changes when
three vertices are kept fixed and the fourth one is moving. For instance,
it was shown in \cite{dv} that the function decreases when we move the
fourth vertex into certain directions.

\subsection{Trapezoid (Burnside \cite{Bur})}
In \cite[pp. 237-239]{bsv} so called square frame, the domain between two concentric squares
with parallel sides, was considered. Such a domain can be split into $8$ similar quadrilaterals, and
we shall study here one such quadrilateral with vertices
$1+hi$, $(h-1)i$, $0$, and $1$, $h>1.$
When $h>1$ we have by \cite[pp.\ 103-104]{Bo}, \cite{Bur}
\begin{equation} \label{sqframe}
\M(Q;1+hi, (h-1)i,0,1)\equiv M(h)\equiv\K(r)/\K(r')
\end{equation}
where
$$
r = \bigg(\frac{t_1-t_2}{t_1+t_2}\bigg)^2\, , \quad
t_1 = \mu_{1/2}^{-1}\left(\frac{\pi}{2c}\right) \, , \quad
t_2 = \mu_{1/2}^{-1}\left(\frac{\pi c}{2}\right)\, , \quad
c = 2h - 1 \, .
$$
Therefore, the quadrilateral can be conformally mapped onto the
rectangle $1+iM(h)$, $i M(h)$, $0$, $1$, with the vertices
corresponding
to each other. It is clear that $h-1 \le M(h) \le h \,.$
The formula (\ref{sqframe}) has the following approximative version
$$
M(h)= h + c +O(e^{-\pi h}),\qquad c=-1/2- \log
2/\pi\approx -0.720636\,,
$$
given in \cite{ps}. As far as we know there is neither
an explicit nor asymptotic formula
for the case when the angle $\pi/4$ of the trapezoid is
replaced by an angle equal to $\alpha \in (0, \pi/2) \,.$

\subsection{Numerical computation of elliptic integrals}
The computation of the elliptic integrals is efficiently carried out by classical
methods available in most programming environments. Numerical estimates for
$\K(r)$, and hence for $\mu_{1/2}(r)$, are obtained very efficiently by the
following recursive method. For $r\in (0,1)$ let
$$
\left\{\begin{array}{ll}
a_0 = 1,&  b_0= r' = \sqrt{1-r^2},\\
a_{n+1} = (a_n+b_n)/2,& b_{n+1} = \sqrt{a_n b_n},\\
\end{array}\right.
$$
Then the sequences $(a_n)$ and $(b_n)$ have the common limit $\pi /(2 \K(r))$,
and, for each $y\in (0,\infty)$ we can approximate $\mu_{1/2}^{-1}(y)$ numerically by
the Newton-Raphson iteration. For details. see e.g. \cite[3.22, 5.32]{avv}
and \cite[2.11]{hvv}.

\section{Validation of algorithms: convex quadrilaterals}

Validation of the algorithms for the modulus of a quadrilateral
will be discussed in two main cases: convex quadrilaterals
and the case of a general polygonal quadrilateral. In this
section the case of a convex quadrilateral will be discussed
for the following three algorithms: (a) the SC Toolbox in MATLAB
written by Driscoll \cite{Dr}, (b) the AFEM software due
to Samuelsson \cite{bsv}, (c) the $hp$-method of the present
paper implemented in the Mathematica language using the double precision. The reference
computation is carried out by the algorithm of \cite{hvv}, implemented
in \cite{hvv} in the Mathematica language 
(the algorithm {\tt QM[A,B]} implementing the formula in Theorem
\ref{modquad}). This implementation makes use of multiple precision
arithmetic for root finding of a transcendental equation
involving the hypergeometric function. All the SC Toolbox tests in this paper were carried out with the settings {\tt  precision =  1e-14}.

\subsection{Setup of the validation test}
All our tests were carried out in the same fashion using the
reciprocal identity (\ref{myrecidty}) and considering a
quadrilateral with the vertices $a,b,0,1$ with $\im a>0$, $\im b>0$, and the line segments
joining the vertices as the boundary arcs. The vertices $b,0,1$
were kept fixed and the vertex $a$ varied over a rectangular
region in the complex plane. The numerical value $b=-0.2+i\cdot1.2$
was used and the lower left (upper right) corner of the rectangular
region was $0.5+i\cdot 0.2$ ($1.5+i\cdot1.2$). Examples of such quadrilaterals, along with some minimal meshes used in the computation, are illustrated in Figure \ref{generatedmeshesfig}. The test functional, based on the
reciprocal identity (\ref{myrecidty}),  is
\begin{equation}
\mathrm{test}(a,b) =\big| \QM(a,b,0,1)\QM\big((b-1)/(a-1),1/(1-a),0,1\big) -1\big|
\end{equation}
which vanishes identically. The values of this test functional are reported in Table \ref{table1} for the fixed value  $b=-0.2+i \cdot 1.2$ when $a$ runs through the aforementioned rectangular region. A table of values of $\QM(m+in,i,0,1)$, $m,n=1,\ldots,5$ is given in \cite[Table~1]{hvv}.

\begin{table}[ht]
\caption{Tests related to the convex quadrilateral, with $(0.15,18)$-meshes used in the $hp$-method.
}\label{table1}
\renewcommand\arraystretch{1}
\noindent

\begin{center}\footnotesize
\begin{tabular}{|l|l|l|}
\hline
Method               & \multicolumn{2}{|c|}{Error range ($\mathrm{test}(a,b)$)} \\
\hline
AFEM                 & $1.15\cdot 10^{-11}$ & $1.41\cdot 10^{-8}$\\
SC Toolbox           & $1.11\cdot 10^{-16}$ & $1.55\cdot 10^{-14}$\\
$hp$-method ($p=12$) & $6.51\cdot 10^{-14}$ & $7.84\cdot 10^{-9}$\\
$hp$-method ($p=15$) & $2.22\cdot 10^{-16}$ & $1.42\cdot 10^{-10}$\\
$hp$-method ($p=18$) & $1.11\cdot 10^{-16}$ & $3.90\cdot 10^{-12}$\\
\hline
\end{tabular}
\end{center}
\end{table}

\subsection{The reference computation}
We used the Mathematica script of \cite{hvv} for solving the
equation in Theorem \ref{modquad} for the computation of $\QM(a,b,0,1)$
in order to carry out the test. The conclusion was that the amplitude of the error was roughly $10^{-17}$ i.e. there was practically no error. Note that the quadrilateral here is not always
convex. On the basis of numerical experiments, it seems that the reference method of \cite{hvv} does also work in non-convex cases, but this has not been rigorously proved.

\begin{figure}
\begin{center}
\includegraphics[width=.31\textwidth]{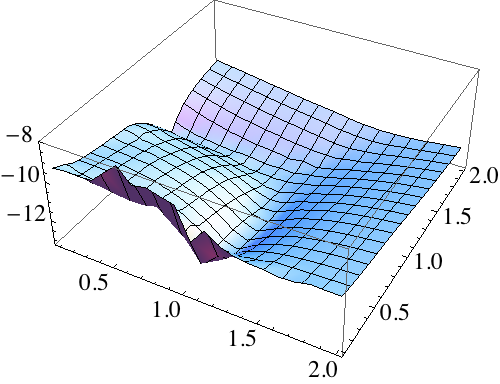}
\includegraphics[width=.31\textwidth]{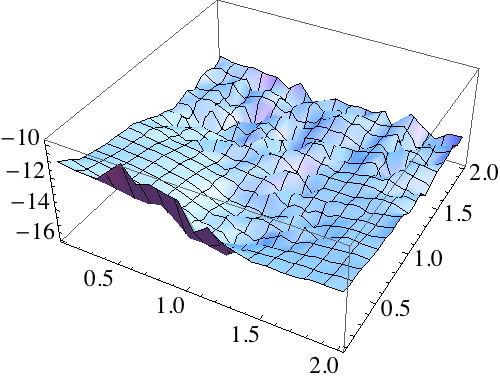}
\includegraphics[width=.31\textwidth]{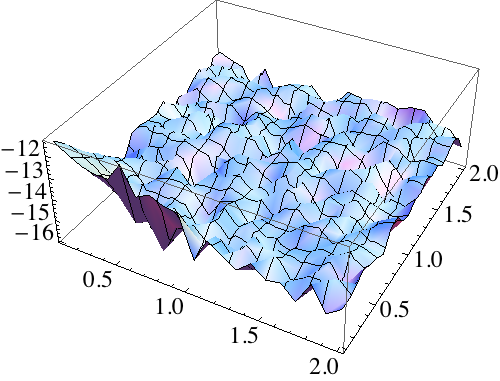}
\end{center}
\caption{Logarithm (with base $10$) of errors over the domain $[0.1,2]\times [0.1,2]$,
corresponding to values of $p=12,15,18$ starting from above. The error estimate is obtained
by using the identity (\ref{myrecidty}).
}\label{errorfigures}
\end{figure}

\section{Validation: polygonal quadrilaterals}

In this section we will consider the validation of the algorithms
for the modulus of a quadrilateral in the case of polygonal domains
with $q>4$ vertices. In the case considered in the previous section
there was a reference computational method, providing the reference
value for the moduli. There is no similar
formula available for the general polygonal case.

\subsection{Setup of the validation test}
All our tests were carried out in the same fashion as in the previous
section, using the reciprocal identity (\ref{myrecidty}). We selected
a quadruple of points $\{z_1,z_2, z_3,z_4\}\,,$ which is a subset of
the set of vertices defining the polygon $D\,,$ and assume that these
are positively oriented. Thus $(D;z_1,z_2, z_3,z_4)$ is a quadrilateral
to which the reciprocal identity (\ref{myrecidty}) applies.

\subsection{The notation ${\rm cmodu}(w,k_1,k_2)$ and ${\rm modu}(w,k_1,k_2)$}
Suppose that $w$ is a vector of $p$ complex numbers such that the
points $w_1, \ldots, w_q$, $q\ge 5,$ are the vertices of a polygon $D$ and that
they define a positive orientation of the boundary. Choose indices
$k_1,k_2 \in \{1,\ldots,p-1\}$ with $k_1<k_2\,$ and set
$z_1=w_{k_1}$, $z_2=w_{k_1+1}$, $z_3=w_{k_2}$, $z_4=w_{k_2+1}\,.$ Then
we define
$$
\mathrm{cmodu}(w,k_1,k_2)= \M(D; z_1,z_2, z_3,z_4)\,, \quad
\mathrm{modu}(w,k_1,k_2)= \M(D; z_2, z_3,z_4, z_1)\, .
$$
By the reciprocal relation (\ref{recipridty}) we have
\begin{equation} \label{myidty}
\mathrm{cmodu}(w,k_1,k_2)\cdot\mathrm{modu}(w,k_1,k_2)=1 \,.
\end{equation}

\subsection{L-shaped region}\label{lshaped}
\begin{figure}[ht]
\begin{center}
\includegraphics[width=.35\textwidth]{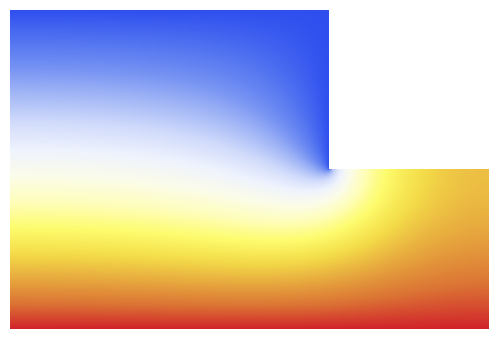}
\includegraphics[width=.35\textwidth]{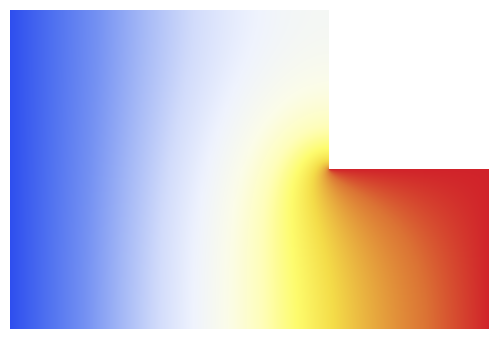}
\end{center}
\caption{Potential functions in the case of L-shaped region \ref{lshaped}. The vertices of the region $Q$ are $z_1=(0,0)$, $z_2=(3,0)$, $z_3=(3,1)$, $z_4=(2,1)$, $z_5=(2,2)$ and $z_6=(0,2)$. Potential functions related to $\M(Q; z_2,z_4,z_6,z_1)\approx 1.5081540958548603$ (left), and $\M(Q; z_1,z_2,z_4,z_6)\approx 0.6630622181450123$ (right), are illustrated.
}\label{lfig}
\end{figure}

\begin{figure}
\begin{center}
\includegraphics[width=.33\textwidth]{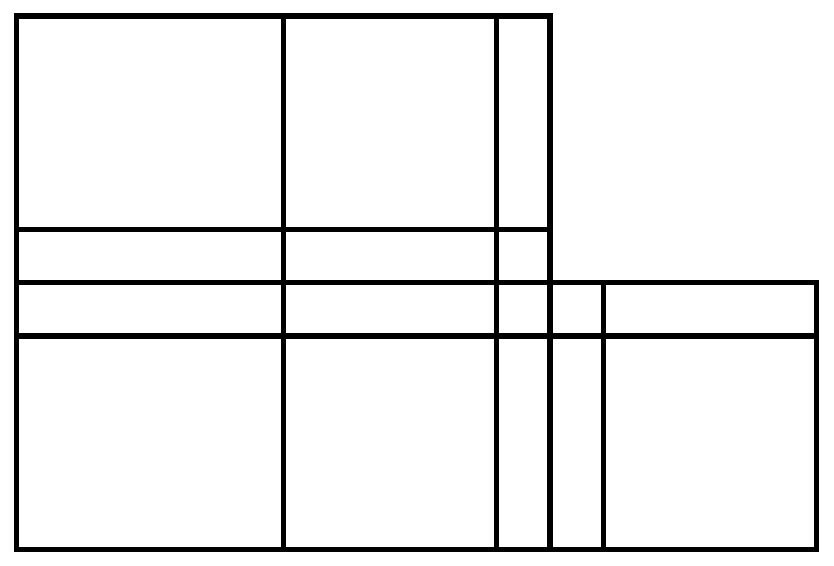}
\includegraphics[width=.33\textwidth]{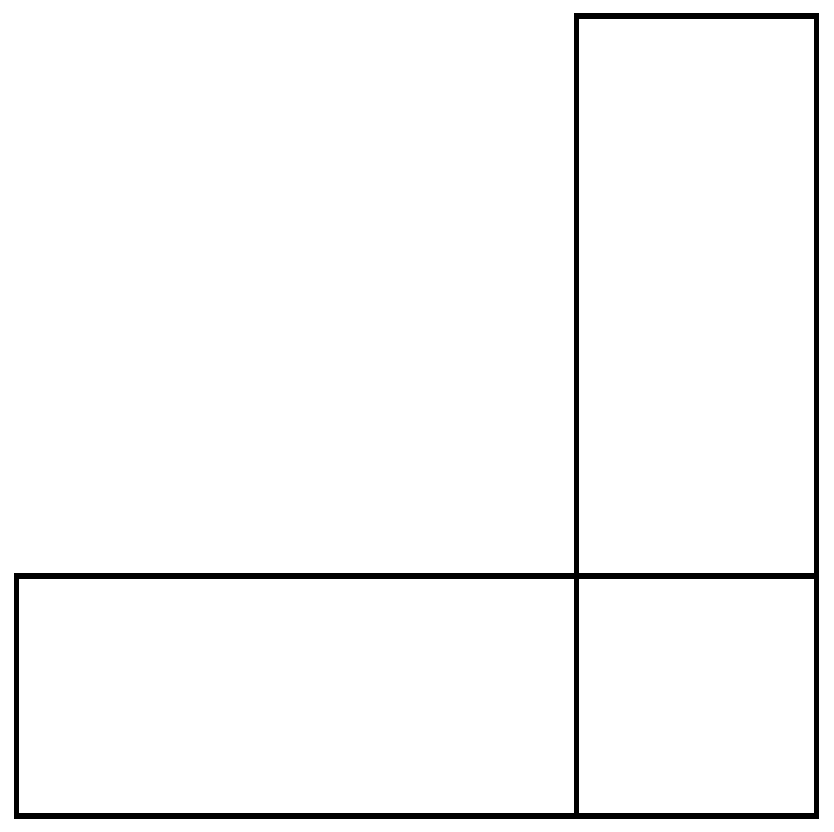}
\end{center}
\caption{Minimal (see \ref{meshing}) meshes for domains of \ref{lshaped} and \ref{square}.
}\label{minmesh}
\end{figure}

\begin{figure}
\begin{center}
\includegraphics[width=.6\textwidth]{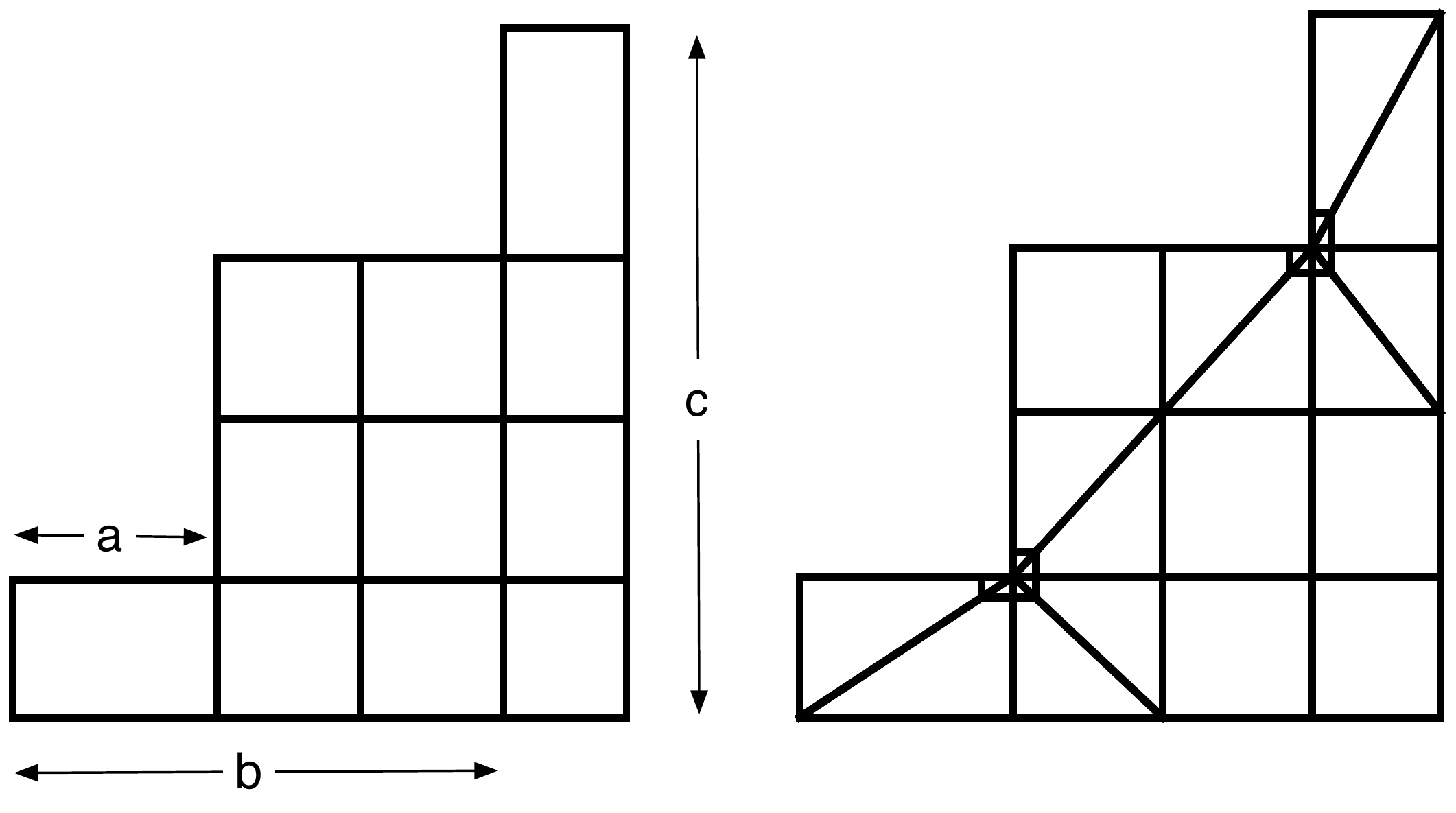}
\end{center}
\caption{Meshing setup for cross in square.}\label{crossmeshfig}
\end{figure}

The L-shaped region:
$$
L(a,b,c,d) =L_1\cup L_2,\quad L_1=\{ z\in \C : 0 < \re(z) < a,\, 0<\im(z) <b \},
$$
$$
L_2= \{ z\in {\mathbb C} : 0 < \re(z) < d,\, 0<\im(z) <c \}\, ,\; 0<d<a, \;0<b<c\,,
$$
is a standard domain considered by several authors for various
computational tasks. In the context of computation of the moduli
it was investigated by Gaier \cite{gai} and we will compare
our results to his results. In the test cases all the vertices had
integer coordinates in the range $[1,4 ]\,.$ Since we consider an integer coordinate domain, simple quadrilateral grid has the desired properties of the minimal mesh, see Figure \ref{minmesh}.
An example of such a mesh is shown in Figure \ref{crossmeshfig}. The results are summarized in Table \ref{table8}, and the potential functions are illustrated by Figure \ref{lfig}.

\begin{table}[ht]
\caption{Tests of (\ref{myidty}) for L-shaped regions (see \ref{lshaped} and Figure \ref{lfig}), with $(0.15,12)$-meshes used in the $hp$-method.}\label{table8}

\begin{center}\footnotesize
\begin{tabular}{|l|l|l|}
\hline
Method               & \multicolumn{2}{|c|}{Error range} \\
\hline
AFEM                 & $1.80\cdot 10^{-10}$ & $7.10\cdot 10^{-10}$\\
SC Toolbox           & $2.22\cdot 10^{-16}$ & $ 2.58\cdot 10^{-14}$\\
$hp$-method ($p=12$) & $4.01\cdot 10^{-11}$ & $1.59\cdot 10^{-10}$\\
$hp$-method ($p=16$) & $8.03\cdot 10^{-13}$ & $2.28\cdot 10^{-12}$\\
$hp$-method ($p=20$) & $5.98\cdot 10^{-13}$ & $1.80\cdot 10^{-12}$\\
\hline
\end{tabular}
\end{center}
\end{table}

\begin{table}[ht]
\caption{Table of capacity values for square in square \ref{square}, $p=16$. Error is the difference to the exact value $4\pi/\mu_{1/2}(r)$, where $r$ is as in (\ref{starfm}).}\label{table3}

\begin{center}\footnotesize
\[
\begin{array}{|l|l|l|l|}
\hline
{{\bf a}}&{{\bf \hbox{Capacity}}}&{{\bf \hbox{Error ($hp$)}}} & {{\bf \hbox{Error (SC)}}}\\
\hline
0.1 & 2.83977741905223   & 2.35 \cdot 10^{-15} & 3.17 \cdot 10^{-14}\\
0.2 & 4.134487024234081  & 1.93 \cdot 10^{-15} & 2.10 \cdot 10^{-15}\\
0.3 & 5.632828000941654  & 1.58 \cdot 10^{-16} & 2.69 \cdot 10^{-16}\\
0.4 & 7.5615315398105745 & 1.17 \cdot 10^{-15} & 1.50 \cdot 10^{-16}\\
0.5 & 10.23409256936805  & 1.74 \cdot 10^{-16} & 3.42 \cdot 10^{-16}\\
0.6 & 14.234879675824363 & 7.49 \cdot 10^{-16} & 1.35 \cdot 10^{-15}\\
0.7 & 20.901581676413954 & 0.                  & 2.89 \cdot 10^{-16}\\
0.8 & 34.23491519877346  & 6.23 \cdot 10^{-16} & 3.61 \cdot 10^{-16}\\
0.9 & 74.23491519877882  & 3.83 \cdot 10^{-16} & 8.31 \cdot 10^{-15}\\
\hline
\end{array}
\]
\end{center}
\end{table}

\begin{table}[ht]
\caption{Table of capacity values for cross in square \ref{crosssquare}, $p=16$. The numerical values and their respective differences are given.}\label{table4}

\begin{center}\footnotesize
\[
\begin{array}{|l|l|l|l|l|l|}
\hline
{{\bf a}}&{{\bf b}}&{{\bf c}}&
{{\bf \hbox{Capacity (SC)}}}& {{\bf \hbox{Capacity ($hp$)}}} &{{\bf \hbox{Difference}}}\\
\hline
0.5&1.2&1.5&
21.94721953515564 &
21.94721953515577 & 5.99 \cdot 10^{-15}\\
0.5&1.0&1.5&
14.00279904484107 &
14.00279904484109 & 8.88\cdot 10^{-16}\\
0.2&0.7&1.2&
9.186926595881523 &
9.186926595881525 & 1.93\cdot 10^{-16}\\
0.1&0.8&1.1&
11.256582318490887 &
11.256582318490889 & 1.58 \cdot 10^{-16}\\
0.5&0.6&1.5&
7.323269585560689 &
7.323269585567927 & 9.88\cdot 10^{-13}\\
0.1&1.2&1.3&
23.13861453810508 &
23.13861453810529  & 8.91  \cdot 10^{-15}\\
\hline
\end{array}
\]
\end{center}
\end{table}

\section{Ring domains}

In this section, we compare $hp$-FEM with exact values and with AFEM and SC Toolbox in certain ring domains. The square in square and cross in square cases were previously considered in \cite{bsv} and numerical values were reported in \cite[Table1, Table 4]{bsv}. Our numerical
results in Tables \ref{table3} and \ref{table4} provide $12$ decimal places whereas in \cite{bsv} only $6$ decimal places were given. Due to the symmetry of the situation it is possible to reduce the computational load for some domains. Of these we discuss here two cases: (a) square in square \ref{square} and (b) cross in square \ref{crosssquare}. These ring domains are symmetric with respect to both the $x$- and $y$-axes, and they are divided into four similar parts by the coordinate axes.

\begin{figure}
\begin{center}
\includegraphics[width=.27\textwidth]{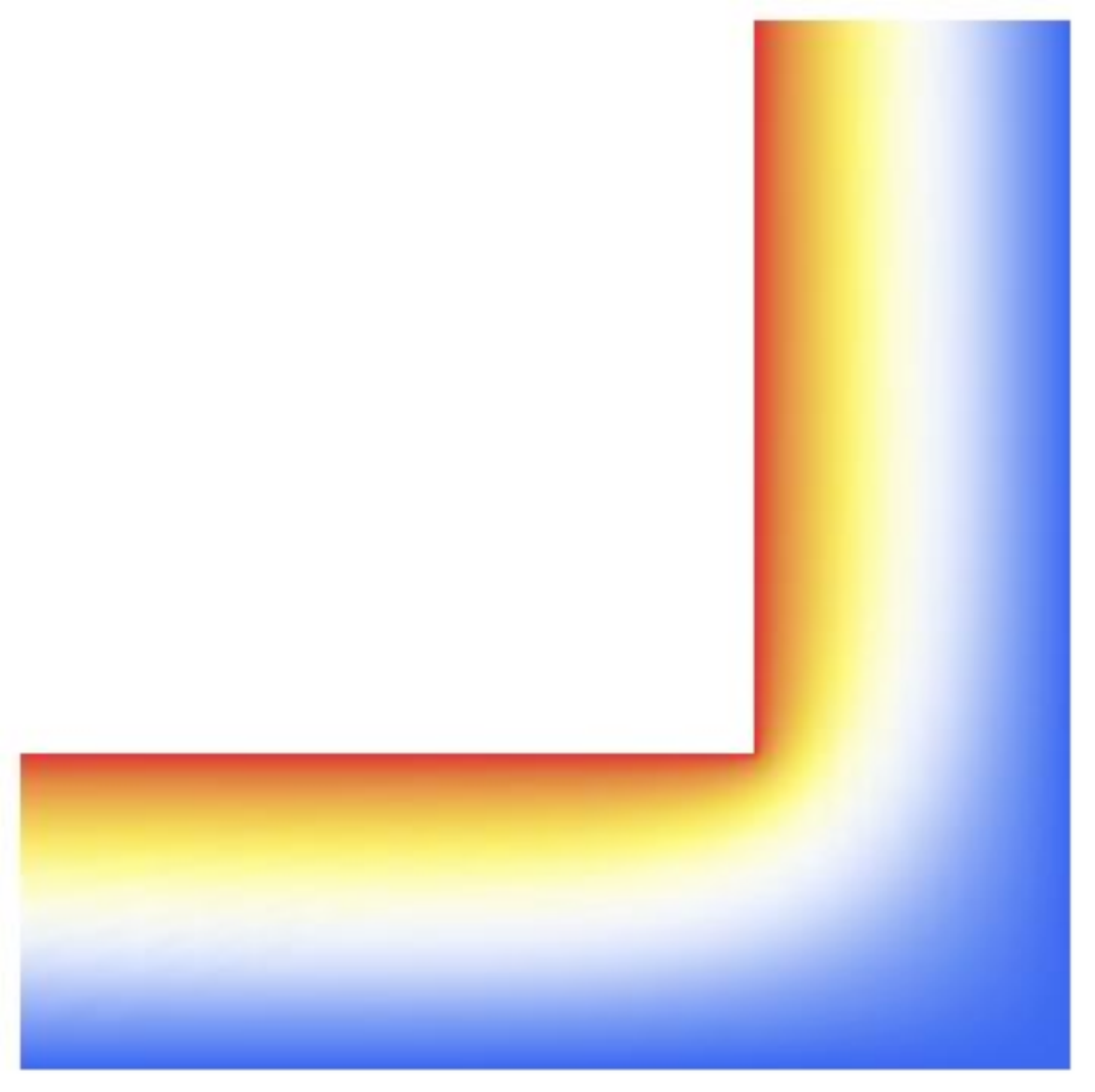}
\includegraphics[width=.27\textwidth]{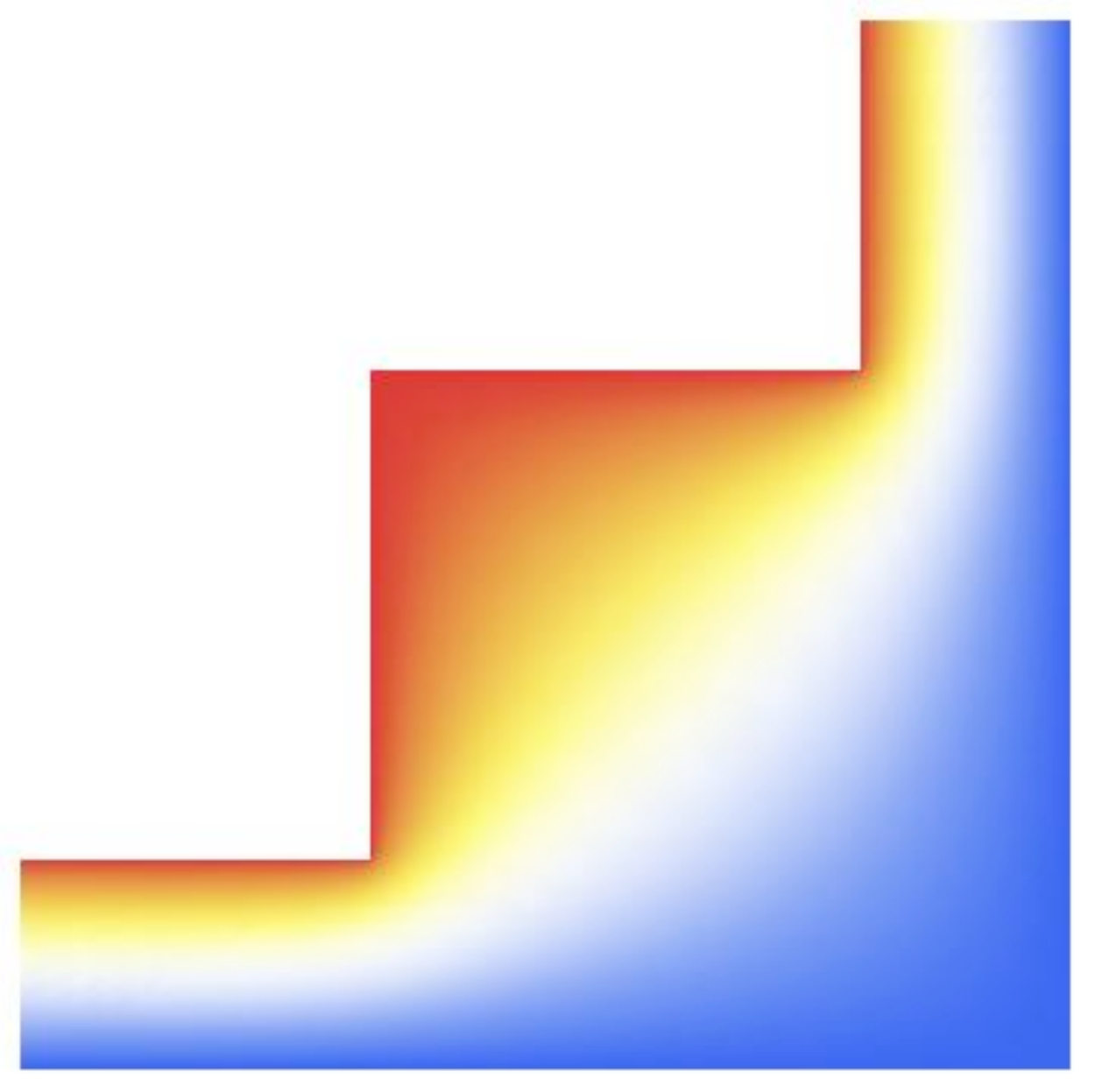}
\end{center}
\caption{Potential functions: square in square and cross in square. Because of the symmetry, only one fourth of the picture is shown.}\label{potfig}
\end{figure}

\subsection{Square in square} \label{square}
We compute here the capacity of the ring domain with plates
$E=[-a,a]\times [-a,a]$ and $F={\C_\infty} \setminus ( (-1,1) \times
(-1,1))$, $0<a<1$. The results with SC and the $hp$-method with $(0.15,16)$-meshes are summarized in
Table \ref{table1}. For computation of the capacity, the ring domain is first split
into four similar quadrilaterals. For the potential function, see Figure \ref{potfig}. Note that
in this case, the exact values of the potential are known, see (\ref{sqframe}) and the
related trapezoid type quadrilateral example. Explicitly, with $c=(1-a)/(1+a)$ and
\begin{equation}
\label{starfm}
  u=  \mu_{1/2}^{-1}\bigg(\frac{\pi \, c}{2}\bigg),\;\; v= \mu_{1/2}^{-1}\bigg(\frac{\pi}{2c}\bigg),\;\;  r= \bigg( \frac{u-v}{u+v}\bigg)^2,
\end{equation}
the capacity equals $4\pi/ \mu_{1/2}(r) \,.$

\subsection{Cross in square} \label{crosssquare} Let $G_{ab}=\{(x,y):|x|\leq a, |y|\leq b\}\cup \{(x,y):|x|\leq b,
|y|\leq a\}$. and $G_c=\{(x,y):|x|<c,|y|<c\}$, where $a<c$ and $b<c$.
We compute the capacity of the ring domain $R=G_c\setminus G_{ab}$.
The results with SC and the $hp$-method with $(0.15,16)$-meshes are summarized in Table \ref{table4}.
For computation of the capacity, the ring domain is again first split into four
similar quadrilaterals. The mesh for the quadrilaterals is given in
Figure \ref{crossmeshfig}, and the potential function is given in Figure \ref{potfig}.
The exact values are not known in this case.

Since the underlying mesh topology
remains constant in both examples above we have computed the
results using exactly the same mesh  template for every subproblem, e.g.
Figure \ref{crossmeshfig} for Cross in square, $a=0.5,b=1.2,c=1.5$,
except for the extremal cases in terms of element distortion
$a=0.9$ for the square in square, and
the case $a = 0.5, b = 0.6, c = 1.5$ for cross in square.
Thus, the results also measure the robustness of the method with respect to
moderate element distortion. Also, in both cases due to symmetry we have graded the mesh only
to the reentrant corners of the domain.

\begin{figure}[ht]
\begin{center}
\includegraphics[width=.40\textwidth]{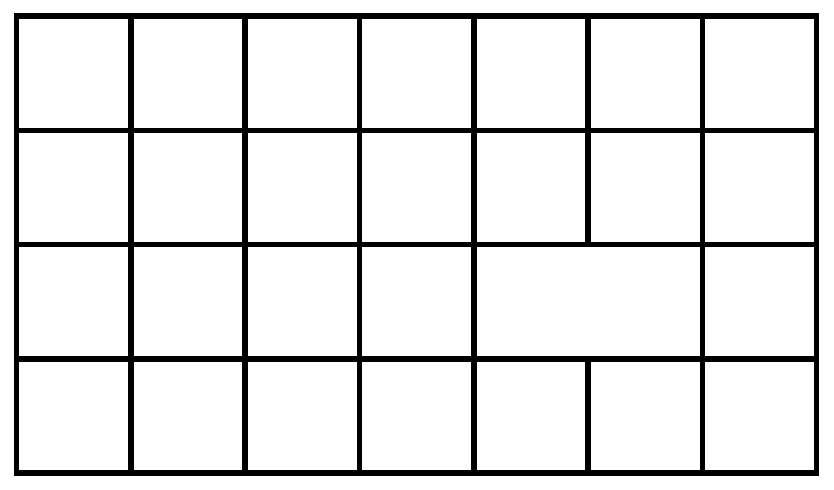}
\includegraphics[width=.40\textwidth]{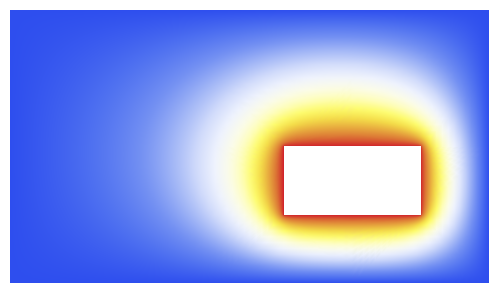}
\end{center}
\caption{Rectangle in rectangle; $\{a,b,c,d\}=\{4,1,6,2\}$, the minimal mesh and the potential function.}\label{rrpotfig}
\end{figure}

\subsection{Rectangle in rectangle}\label{rectrect}
 Let $G_{abcd}=\{(x,y):a\leq x \leq c, b\leq y \leq d\}$ and $G=\{(x,y):0\leq x\leq7,\ 0\leq y\leq4\}$.
We compute the capacity of the ring domain $R=G\setminus G_{abcd}$.
Here we consider a subset of possible cases when $a,b,c,d \in \mathbb{N}$.
The results computed using
the $hp$-method with $(0.15,16)$-meshes are summarized in Table \ref{rrtable}.
The potential function for the case $\{a,b,c,d\}=\{4,1,6,2\}$ is given in Figure \ref{rrpotfig}.
The exact values are not known in this case.

Again, we have employed the same mesh template (simple quadrilateral grid as in Figure \ref{crossmeshfig}) over the entire test set. Grading has been used in the corners of $G_{abcd}$ only. From results of Table \ref{rrtable} we can also see that some of the configurations are
symmetric in terms of capacity. In these cases the differences in the computed
values are less than $10^{-13}$.

\begin{table}[ht]
\caption{Table of capacity values for rectangle in rectangle \ref{rectrect}.}\label{rrtable}
\begin{center}\footnotesize
\begin{tabular}{|l|l|l|l| l | l|}
\hline
$a$ & $b$ & $c$ & $d$ & $p$ & Capacity\\
\hline
1 & 1 & 2 & 2 & 20 & 5.210320385649294 \\
1 & 1 & 3 & 2 & 19 & 6.746053277945276\\
1 & 1 & 4 & 2 & 20 & 8.27007839293125 \\
1 & 1 & 5 & 2 & 19 & 9.86240917550835\\
1 & 1 & 6 & 2 & 17 & 11.89718127369752\\
2 & 1 & 3 & 2 & 18 & 4.692072335693745\\
2 & 1 & 4 & 2 & 18 & 6.232078709256309\\
2 & 1 & 5 & 2 & 20 & 7.827105378062926\\
2 & 1 & 6 & 2 & 17 & 9.86240917550835\\
3 & 1 & 4 & 2 & 17 & 4.621123827863167\\
3 & 1 & 5 & 2 & 20 & 6.232078709256313\\
3 & 1 & 6 & 2 & 18 & 8.2700783929313\\
4 & 1 & 5 & 2 & 19 & 4.69207233569376\\
4 & 1 & 6 & 2 & 20 & 6.746053277945233\\
5 & 1 & 6 & 2 & 20 & 5.210320385649318\\
\hline
\end{tabular}
\end{center}
\end{table}

\begin{figure}
\begin{center}
\includegraphics[width=0.3\textwidth]{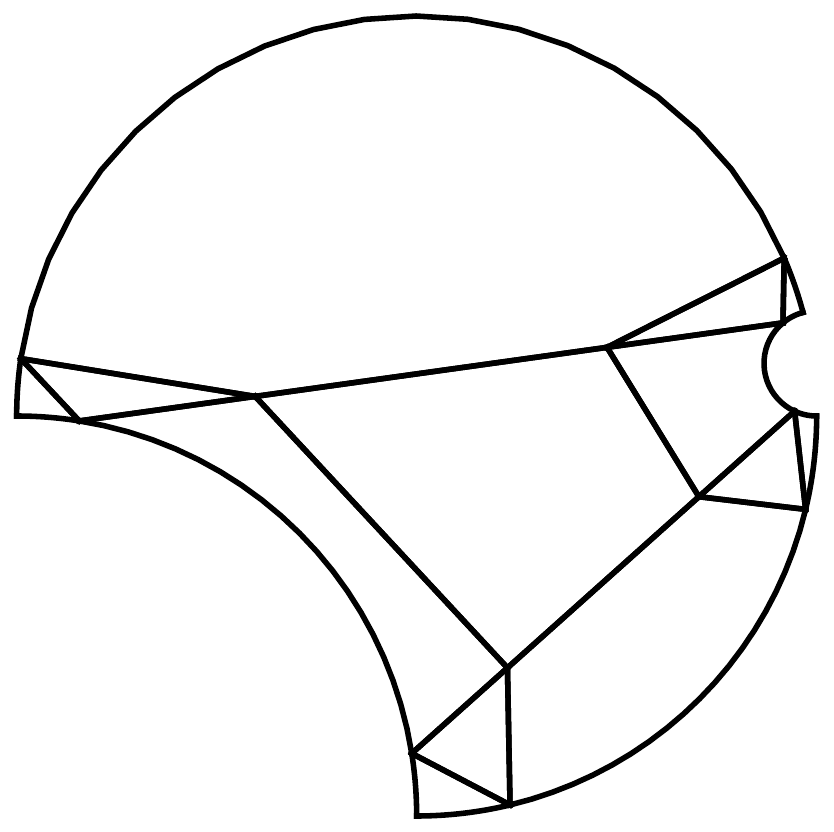}
\includegraphics[width=0.3\textwidth]{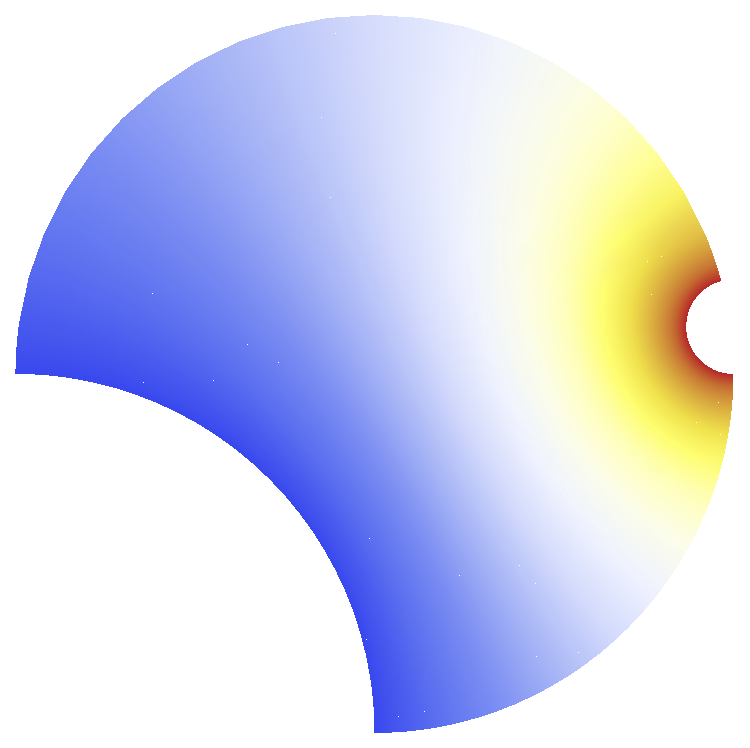}
\end{center}
\caption{The quadrilateral $(Q_A;\pi/12, \pi, 3 \pi/2,1)$, the (final) mesh and the potential function. The relative error is $1.23\cdot 10^{-13}$.}\label{qa1}
\end{figure}

\begin{figure}
\begin{center}
\includegraphics[width=0.3\textwidth]{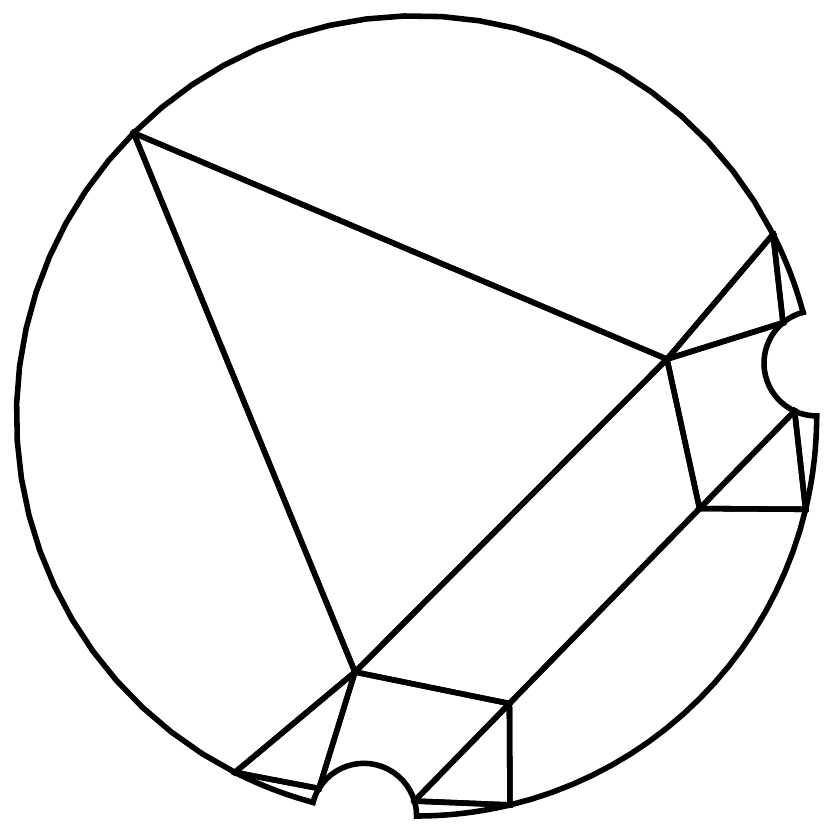}
\includegraphics[width=0.3\textwidth]{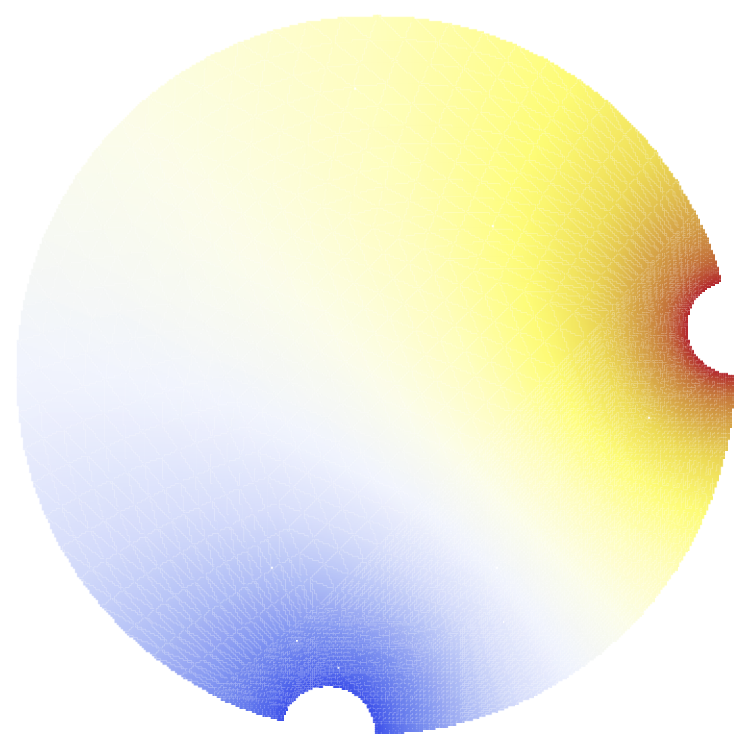}
\end{center}
\caption{The quadrilateral $(Q_A;\pi/12, 17\pi/12, 3 \pi/2,1)$, the (final) mesh and the potential function. The relative error is $6.14\cdot 10^{-14}$.}\label{qa2}
\end{figure}

\begin{figure}
\begin{center}
\includegraphics[width=0.3\textwidth]{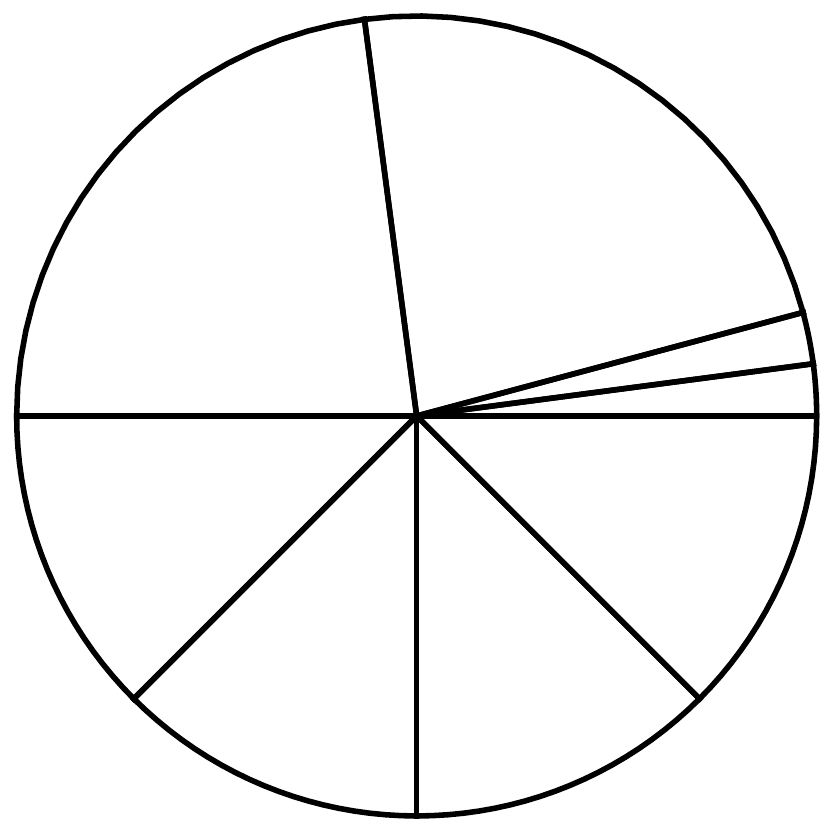}
\includegraphics[width=0.3\textwidth]{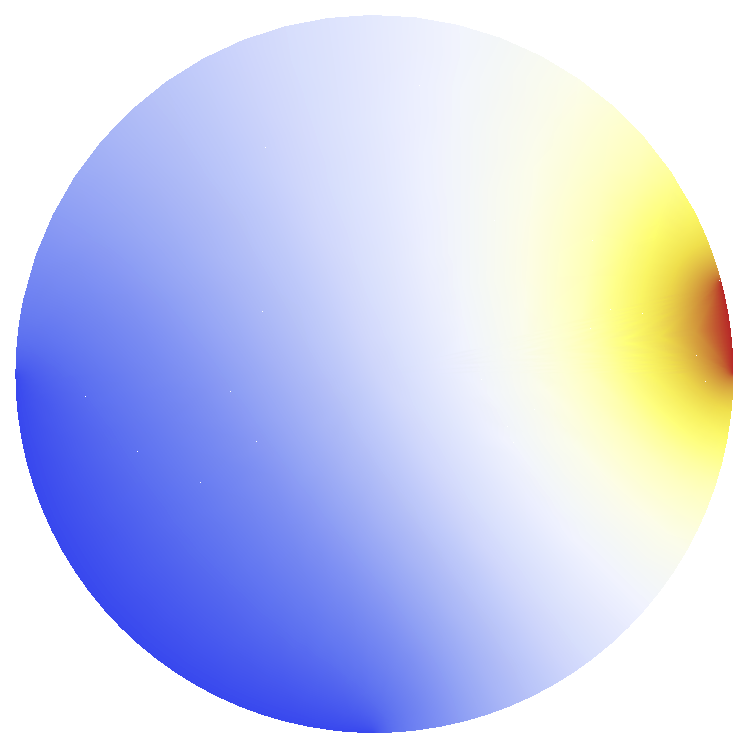}
\end{center}
\caption{The quadrilateral $(Q_B;\pi/12, \pi, 3 \pi/2,1)$, the minimal mesh and the potential function. The relative error is $3.38\cdot 10^{-11}$.}\label{qb1}
\end{figure}

\begin{figure}
\begin{center}
\includegraphics[width=0.3\textwidth]{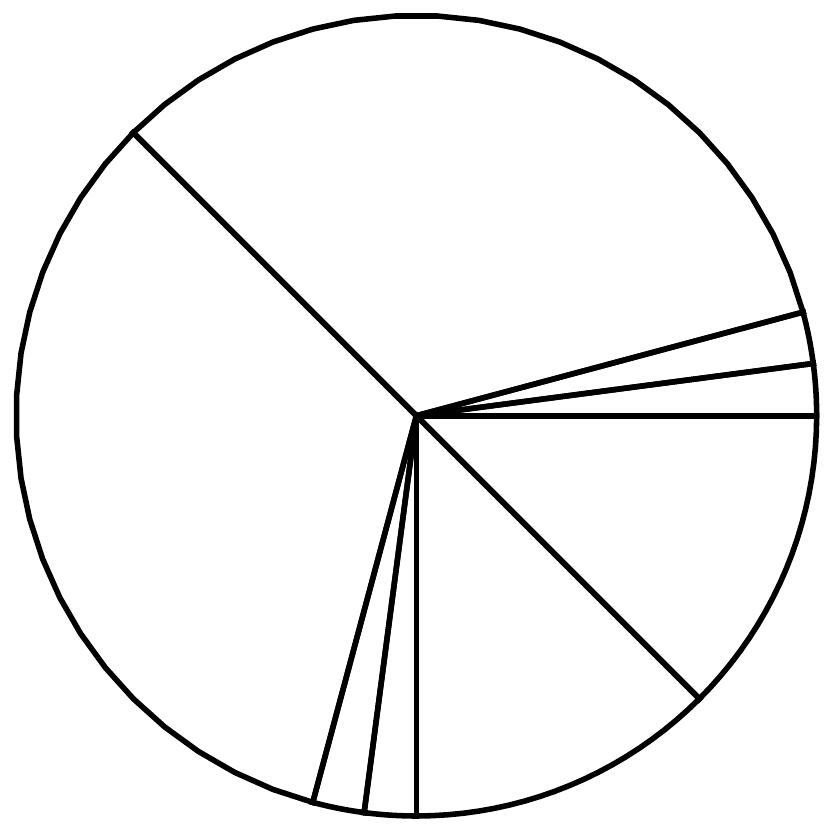}
\includegraphics[width=0.3\textwidth]{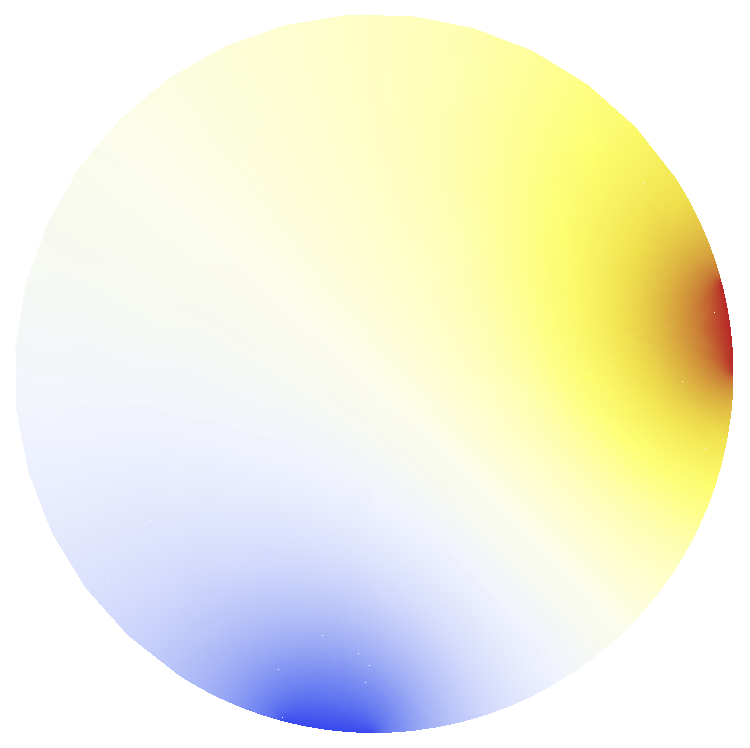}
\end{center}
\caption{The quadrilateral $(Q_B; \pi/12, 17\pi/12, 3 \pi/2,1)$, the minimal mesh and the potential function. The relative error is $5.17\cdot 10^{-11}$.}\label{qb2}
\end{figure}

\begin{figure}
\begin{center}
\includegraphics[width=0.4\textwidth]{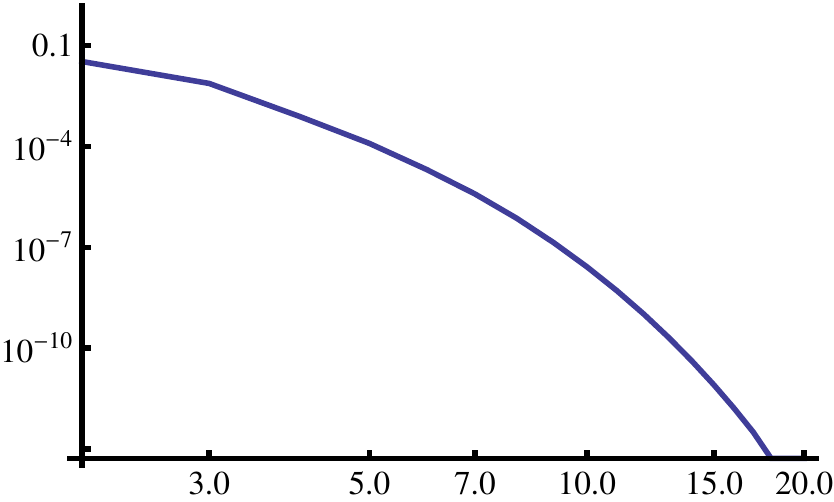}
\includegraphics[width=0.4\textwidth]{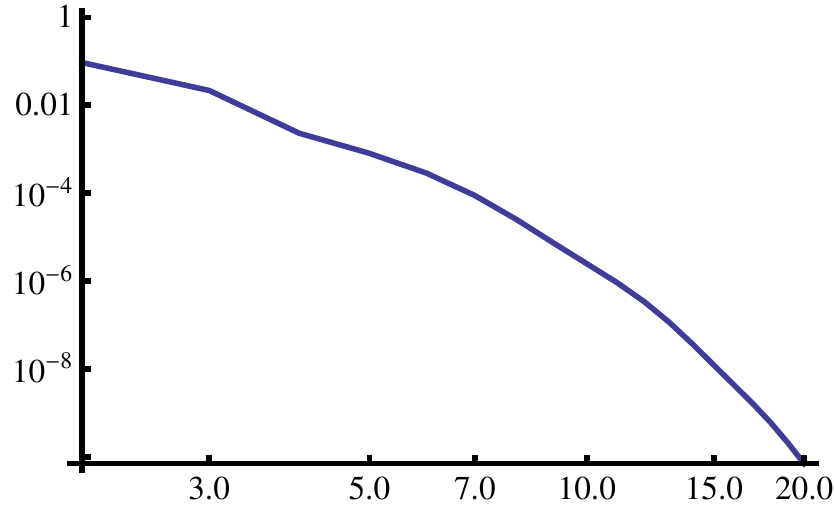}
\end{center}
\caption{The $p$-convergence of the reciprocal error for the quadrilaterals $(Q_A;\pi/12, 17\pi/12, 3 \pi/2,1)$ (left) and
$(Q_B;\pi/12, 17\pi/12, 3 \pi/2,1)$ (right). Logarithmic scale.}\label{qabconvergence}
\end{figure}

\begin{figure}
\begin{center}
\includegraphics[width=0.45\textwidth]{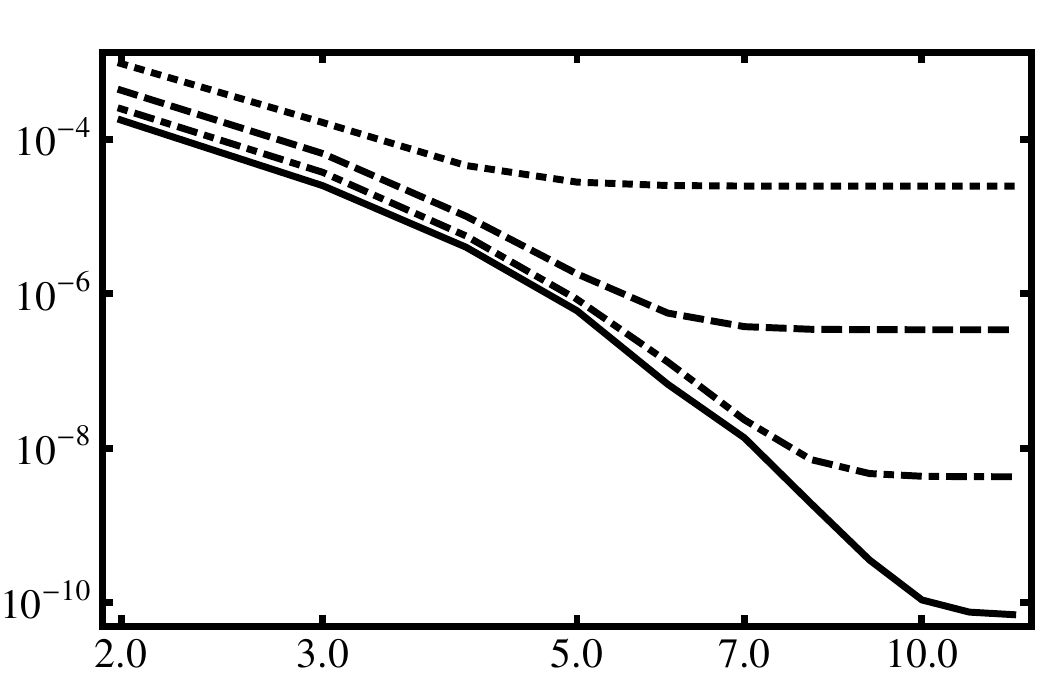}
\end{center}
\caption{The $\nu$-convergence of the reciprocal error for the quadrilateral
$(Q_B;\pi/12, \pi, 3 \pi/2,1)$ for different values of $p$;
$p=8$ dotted line, $p=12$ dashed line, $p=16$ dot-dashed line, $p=20$ solid line,
$\nu=2,\ldots,12$, $\alpha=0.15$. Logarithmic scale.}\label{qnuconvergence}
\end{figure}

\section{Domains with curved boundaries}

In this section, we give further examples featuring domains with curved boundaries. Simple examples of such domains are domains, where four or more points are connected with circular arcs. Some examples related to numerical methods and Schwarz-Christoffel formula for such domains can be found in the literature, e.g. \cite{Br,Hough}. Our method has the advantage that even more general quadrilaterals can be considered, as illustrated by examples given below. Here the meshing has been tuned by monitoring the rate of convergence in the polynomial degree. Both the minimal mesh and the scaling factor have been adjusted until exponential convergence in $p$ has been observed.
Let us consider the quadrilaterals $(Q_A;\pi/12, 17\pi/12, 3 \pi/2,1)$ and
$(Q_B;\pi/12, 17\pi/12, 3 \pi/2,1)$.
In Figures \ref{qabconvergence} and \ref{qnuconvergence} the effects
of choosing the polynomial order and the nesting level, respectively, are shown.
Since the boundary segments of $(Q_A;\pi/12, 17\pi/12, 3 \pi/2,1)$ are
orthogonal, there is no need to refine the mesh. For $(Q_B;\pi/12, 17\pi/12, 3 \pi/2,1)$
the effect of nesting for a given $p$ eventually diminishes as the errors from
outside the corners (sometimes referred to as the modelling error) start to dominate.
On the other hand, the smallest error for $p=20$ is obtained at $\nu=12$.
This error can only be made smaller by modifying the minimal mesh and/or the value of the grading parameter $\alpha$.

\subsection{Circular quadrilaterals}
The absolute ratio of four points $a,b,c,d \in \C$ is defined as
\begin{equation}\label{absratio}
|a,b,c,d| = \frac{|a-c||b-d]}{|a-b||c-d|}.
\end{equation}
The main property of the absolute ratio is the M\"obius invariance:
\begin{equation}\label{mobinv}
|a,b,c,d| = \big|w(a),w(b),w(c),w(d)\big|,
\end{equation}
if $w$ is a M\"obius transformation
\begin{equation}\label{mobtrans}
w(z) = \frac{k z + l}{m z + n},\quad(k n - m l \neq 0).
\end{equation}
Given $z_1, z_2, z_3$ on a circle (or on a line) and $w_1, w_2, w_3$
on a circle (or on a line), there exists a M\"obius transformation $w$ with
$w(z_j)=w_j,\ j=1,2,3$.

\begin{table}[ht]
\caption{Moduli of quadrilaterals $(Q_A;e^{i m \pi/24} ,e^{i n \pi/24}, e^{i r
\pi/24} ,1)$  for several integer triples $(m,n,r)$ computed with the $hp$-method, $p=20$.}\label{table5}

\begin{center}\footnotesize
\begin{tabular}{|l|l|l|l|l|}
\hline
Nodes & Reference & Computed value& Relat. error & Recipr. error \\ \hline
$(2,10,12)$& $0.7071508111121534$ &  $0.7071508111121347$ & $2.64\cdot 10^{-14}$ & $1.02\cdot 10^{-13}$\\
$(2,10,14)$& $0.8074514311467651$ & $0.8074514311467831$ & $2.23\cdot 10^{-14}$ & $2.55\cdot 10^{-14}$\\
$(4,12,18)$& $1.0383251171675787$ & $1.0383251171675796$ & $8.55\cdot 10^{-16}$  & $1.44\cdot 10^{-15}$\\
$(6,16,24)$& $1.170060906774661$ & $1.1700609067746603$ & $5.69\cdot 10^{-16}$  & $2.22\cdot 10^{-15}$\\
$(8,22,32)$& $1.313262425617007$ & $1.3132624256170076$ & $5.07\cdot 10^{-16}$ & $3.33\cdot 10^{-16}$\\
\hline
\end{tabular}
\end{center}
\end{table}

\subsection{Type A}
Let us first consider a quadrilateral whose sides are circular arcs of
intersecting orthogonal circles, i.e.,  angles are $\pi/2$.
Let $0<a<b<c<2\pi$ and choose the points
$\{1, e^{ia}, e^{ib}, e^{ic}\}$ on the unit circle with the absolute ratio
\begin{equation}
\label{udef}
\big|1, e^{ia}, e^{ib}, e^{ic}\big| = \frac{\sin(b/2) \sin((c-a)/2)}{\sin(a/2) \sin((c-b)/2)} = u.
\end{equation}
Let $Q_A$ stand for the domain which is obtained from the unit disk
by cutting away regions bounded by the two orthogonal arcs with endpoints
$\{ 1, e^{ia} \}$ and $\{ e^{ib}, e^{ic} \}\,,$ respectively. Then $Q_A$
determines a quadrilateral
$(Q_A ; e^{ia}, e^{ib}, e^{ic},1)\,.$
Using a suitable M\"obius transformation and
the invariance (\ref{mobinv}) we can map $Q_A$ onto the upper half of the
annulus $\{z \in \C:1<|z|<t \}$ and we obtain the following formula:
\begin{equation}
\M(Q_A ; e^{ia}, e^{ib}, e^{ic}, 1) = \pi / \log t,
\end{equation}
i.e. a half of the modulus of the full annulus, where
\[
t=2u - 1 + 2 \sqrt{u^2-u},\qquad t>1.
\]
The results are summarized in Table \ref{table5}.

\begin{table}[ht]
\caption{Moduli of quadrilaterals $(Q_B;e^{i m \pi/24} ,e^{i n \pi/24}, e^{i r
\pi/24} ,1)$  for several integer triples $(m,n,r)$ computed with the $hp$-method, $p=20$.}\label{table6}

\begin{center}\footnotesize
\begin{tabular}{|l|l|l|l|l|}
\hline
Nodes & Reference & Computed value& Relat. error & Recipr. error \\ \hline
$(2,10,12)$& $0.5389714947317054$ & $0.5389714947624924$ & $5.71\cdot 10^{-11}$ & $7.68\cdot 10^{-11}$\\
$(2,10,14)$& $0.5953434982171909$ & $0.5953434982359955$ & $3.16\cdot 10^{-11}$ & $4.26\cdot 10^{-11}$\\
$(4,12,18)$& $0.7121629047455362$ & $0.7121629047457778$ & $3.39\cdot 10^{-13}$ & $6.06\cdot 10^{-13}$ \\
$(6,16,24)$& $0.7718690862645192$ & $0.7718690862646902$ & $2.22\cdot 10^{-13}$ & $4.09\cdot 10^{-13}$\\
$(8,22,32)$& $0.8319009599091923$ & $0.8319009599093506$ & $1.90\cdot 10^{-13}$ & $3.48\cdot 10^{-13}$\\
\hline
\end{tabular}
\end{center}
\end{table}

\begin{table}[ht]
\caption{Timing for $(Q_B;e^{i\pi/12}, e^{i\pi}, e^{i3 \pi/2},1)$, with $(0.15,12)$-meshes (Apple Mac Pro 2009 Edition 2.26 GHz, Mathematica 7.0.1).} \label{tbl:timing}

\begin{center}\footnotesize
\begin{tabular}{|l|l|l|l|l|}\hline
$p$ & Reciprocal error  &Number of d.o.f. &Time (seconds)  &Setup/solve \\ \hline 
8   &  $1.6 \cdot 10^{-5}$ &  6817 & 4.1 & 9.2 \\
12  &  $2.2 \cdot 10^{-7}$ & 15025 & 9.7 & 3.9 \\
16  &  $2.8 \cdot 10^{-9}$ & 26433 &  27 & 3.4 \\
20  & $4.5 \cdot 10^{-11}$ & 41041 &  67 & 3.3 \\ \hline
\end{tabular}
\end{center}
\end{table}

\subsection{Type B}
Next we let the sides of the quadrilateral be circular arcs be of the unit disk,  and in this case all the angles are equal to $\pi$. Now the unit disk, together with the
boundary points $ e^{ia}, e^{ib}, e^{ic}, 1$ determines a quadrilateral which we denote by
$Q_B\,.$ Using an auxiliary M\"obius transformation of the unit disk onto the
upper half plane we can readily express the modulus using the capacity of the
Teichm\"uller ring domain \cite[Section 7]{avv} and express it as
follows
\begin{equation}
\M(Q_B ;  e^{ia}, e^{ib}, e^{ic},1) = \frac12 \tau(u-1),
\end{equation}
where $u$ is as in (\ref{udef}), and
$$\tau(t) = \pi/\mu_{1/2}(1/\sqrt{1+t})\,,\qquad t >0,\,$$
and $\mu_{1/2}(r)$ is as in (\ref{mu_a}), gives the conformal capacity of the plane Teichm\"uller ring. The results are summarized in Table \ref{table6}.

In every test case, the local stiffness matrices have been integrated, then assembled into the system matrix $A$, and finally two linear systems of equations derived from $A$ have been solved. We present timing results for the $(Q_B;e^{i\pi/12}, e^{i\pi}, e^{i3 \pi/2},1)$ with $\nu = 12$ in Table \ref{tbl:timing}. This case was chosen because, due to the curved geometry, in terms of numerical integration it represents the worst case. The total execution time in seconds and
a dimensionless ratio, system setup time / time spent in the linear solver, are given for the values of $p=8,12,16,20$. In this kind of experiments, specifying the domain and the initial mesh are the most time consuming parts as the execution times are at most minutes and for reasonable accuracy (as in $p=12$) seconds. As one would expect in Mathematica environment, in the range
of problems considered, the system setup time is much longer than time spent in solving the linear systems.

\begin{figure}[ht]
\begin{center}
\includegraphics[width=0.21\textwidth]{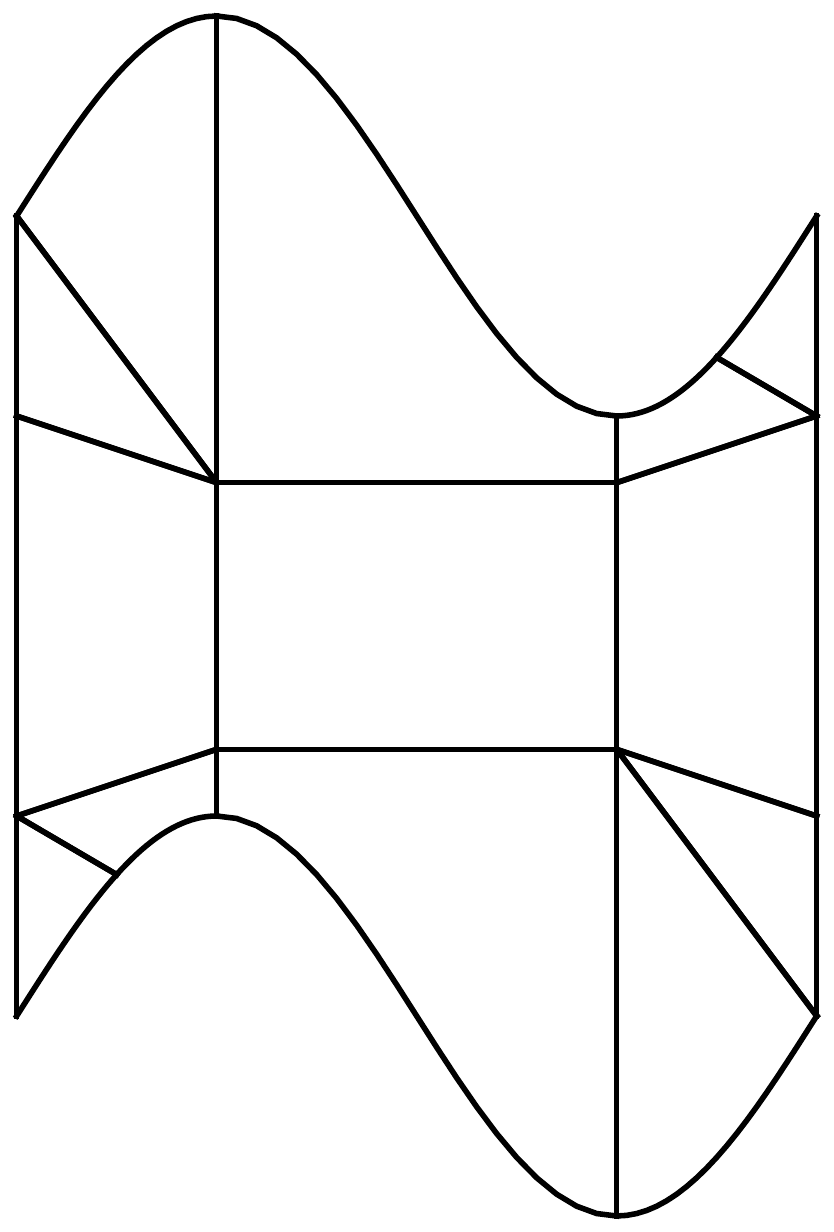}
\includegraphics[width=0.21\textwidth]{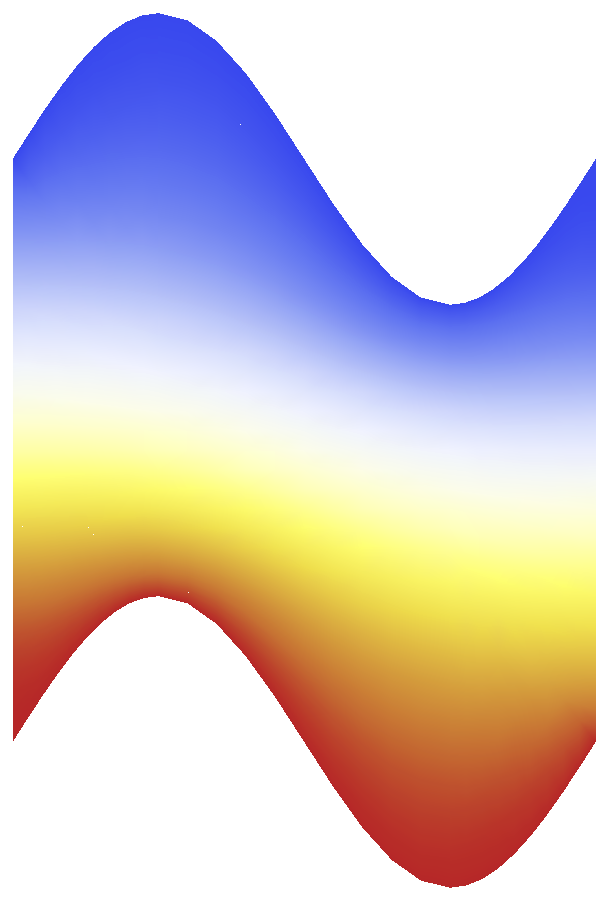}
\end{center}
\caption{Wave: the minimal mesh and the potential function for \ref{wave}.
}\label{wavefig}
\end{figure}

\begin{figure}[ht]
\begin{center}
\includegraphics[width=0.3\textwidth]{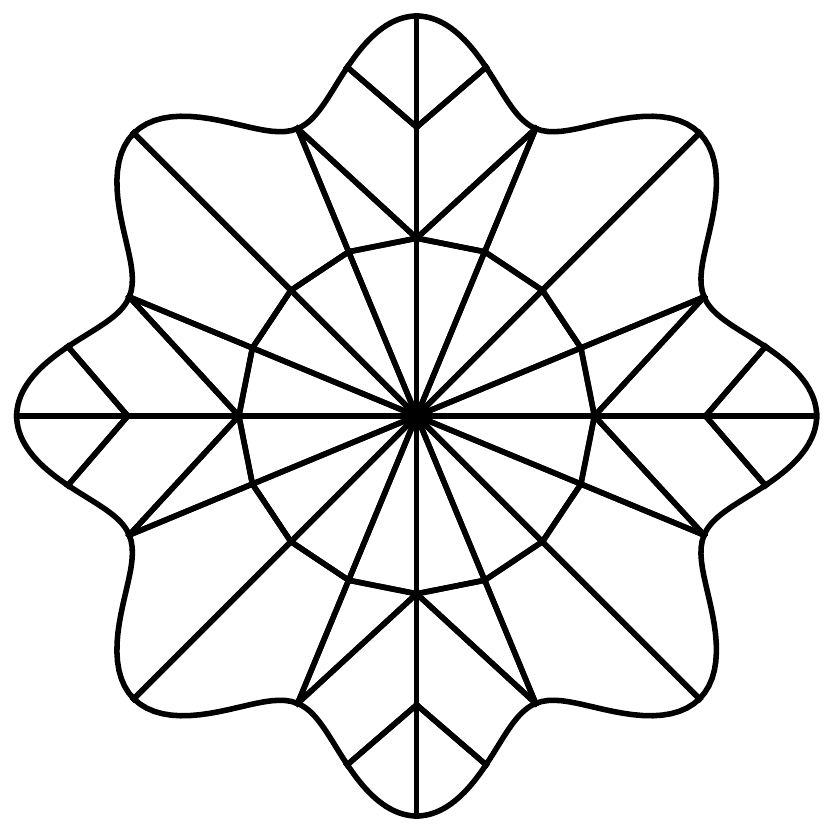}
\includegraphics[width=0.3\textwidth]{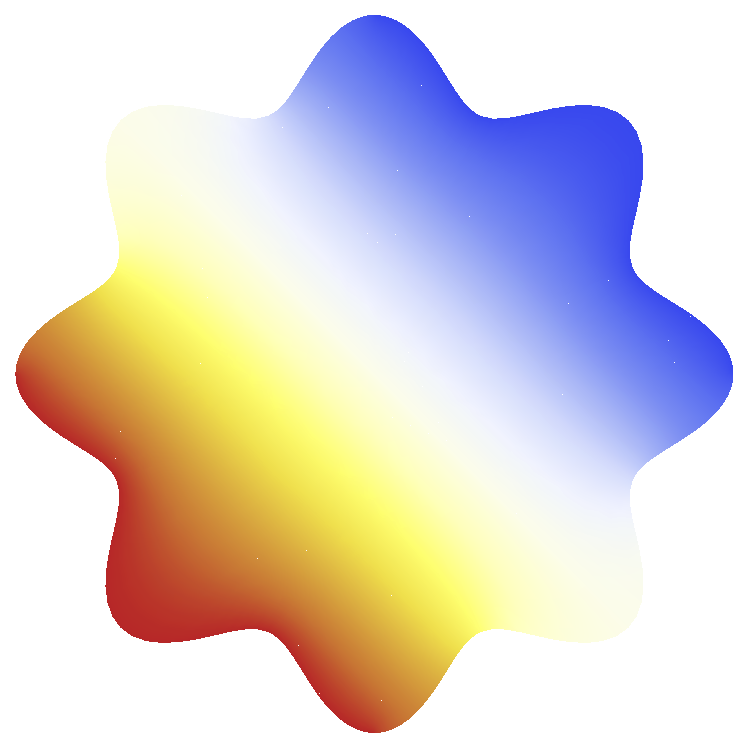}
\end{center}
\caption{Flower I: the minimal mesh and the potential function. See also \ref{flowers} and Table \ref{flowertaba}.
}\label{flowera}
\end{figure}

\begin{figure}[ht]
\begin{center}
\includegraphics[width=0.3\textwidth]{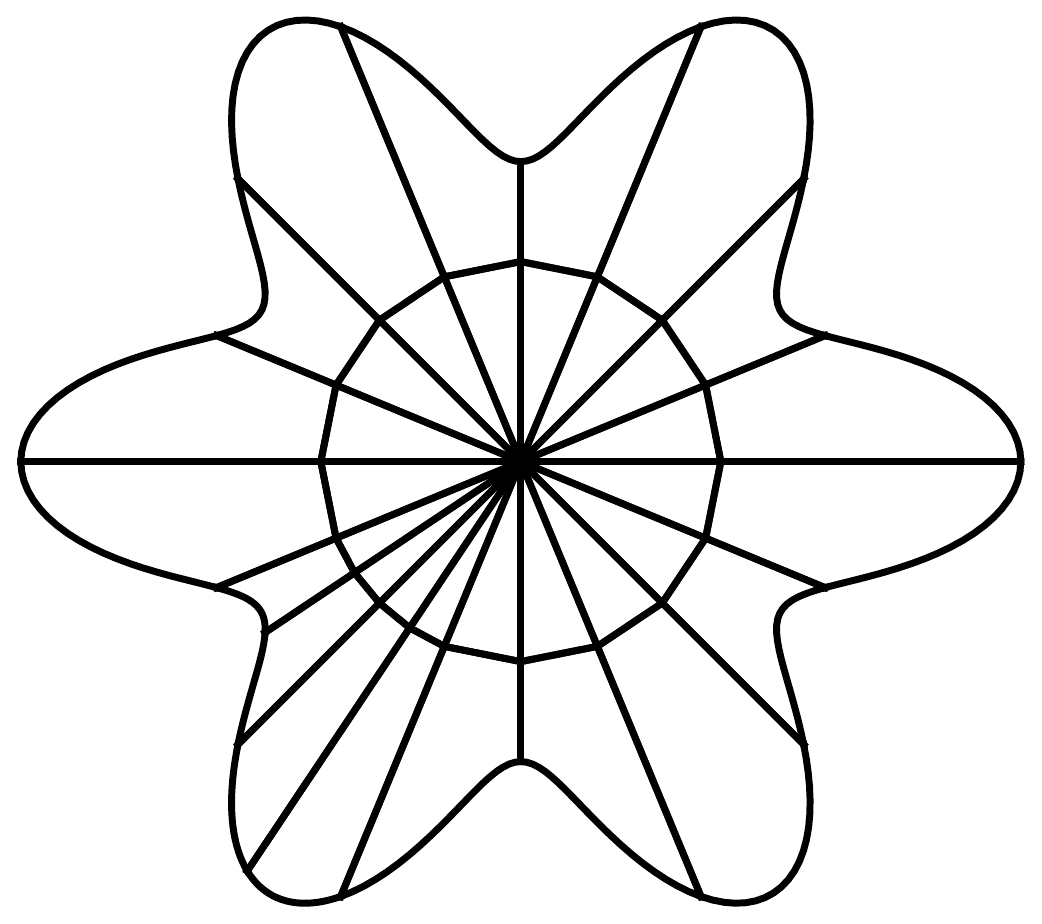}
\includegraphics[width=0.3\textwidth]{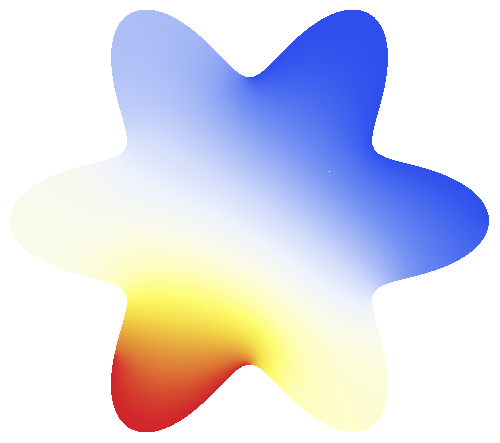}
\end{center}
\caption{Flower II: the minimal mesh and the potential function. See also \ref{flowers} and Table \ref{flowertabb}.
}\label{flowerb}
\end{figure}

Next we consider non-convex examples featuring quadrilaterals with curved boundaries which are not circular segments.

\subsection{Wave}\label{wave}
Let $Q= \{  (x,y) : 0 < x <1,   \sin(2\pi x)/4 <y  < 1+  \sin(2\pi x)/4 \}$,
and $z_1=(0,0)$, $z_2= (1,0)$, $z_3=(1,1)$, $z_4=(0,1).$ Then the $hp$-method with $p=20$ gives
$  \M(Q; z_2, z_3, z_4, z_1)  \approx   1.285385932609546.$ An error estimate based on the reciprocal identity (\ref{recipridty}) is $2.66\cdot 10^{-15}$. For visualization, see Figure \ref{wavefig}.

\subsection{Flowers}\label{flowers}
Let $Q$ be the domain bounded by the curve  $r(\theta)= 0.8 +t \cos(n \pi \theta)$, $0 \le \theta\le 2 \pi$, $n=4,6,8$ and $t=0.1$ or $t=0.2\,.$ Domains of this type are illustrated in Figures \ref{flowera} and \ref{flowerb}. We compute moduli of quadrilaterals $\M(Q; z_1,z_2, z_3, z_4)$, where $z_j=r(\theta_j)$. We consider flower shaped quadrilaterals of type I with $\theta_j=(j-1)\pi/2$ for $j=1,2,3,4$, and type II, where $\theta_1,\theta_2,\theta_3$ are as before, and $\theta_4=5\pi/4$ (see Figures \ref{flowera} and \ref{flowerb}).  The numerical results are summarized in Tables \ref{flowertaba} and \ref{flowertabb}.

\begin{table}[ht]
\caption{Moduli of flower-shaped  quadrilaterals $(t,n)$ of type I computed with the $hp$-method, $p=20$. Note that, because of symmetry, it follows from (\ref{recipridty}) that the exact value of modulus is $1$.
}\label{flowertaba}
\begin{center}\footnotesize
\begin{tabular}{|l|l|l|}
\hline
$n$ & Error ($t=0.1$) & Error ($t=0.2$) \\
\hline
$4$ &  $3.18\cdot 10^{-14}$ & $2.25\cdot 10^{-14}$\\
$6$ &  $3.74\cdot 10^{-11}$ & $8.45\cdot 10^{-11}$\\
$8$ &  $1.34\cdot 10^{-13}$ & $6.27\cdot 10^{-11}$\\
\hline
\end{tabular}
\end{center}
\end{table}

\begin{table}[ht]
\caption{Moduli of flower-shaped  quadrilaterals $(t,n)$ of type II computed with the $hp$-method, $p=20$. The error estimate is obtained by using the reciprocal identity (\ref{recipridty}).
}\label{flowertabb}

\begin{center}\footnotesize
\begin{tabular}{|c|c|c|c|}
\hline
$n$ & $t$ & Error & Modulus \\ \hline
4 & 0.1 & $2.00\cdot 10^{-15}$ & 0.8196442147286799 \\
4 & 0.2 & $1.40\cdot 10^{-13}$ & 0.8196441884805612 \\
6 & 0.1 & $2.34\cdot 10^{-14}$ & 0.7896695654987764\\
6 & 0.2 & $1.43\cdot 10^{-10}$ & 0.7690460663235661\\
8 & 0.1 & $9.05\cdot 10^{-14}$ & 0.8196441884804566\\
8 & 0.2 & $1.38\cdot 10^{-10}$ & 0.8196441885295815\\
\hline
\end{tabular}
\end{center}
\end{table}


\section{Summary}

The computation of the moduli of quadrilaterals and ring domains with piecewise smooth boundaries is a problem frequently occurring in various applications, see \cite{ps2}. There is no general method for such computations
except perhaps the case of polygonal quadrilaterals when the SC Toolbox \cite{Dr,DrTr} may be considered as the ``state-of-the-art'' tool. For the case of ring domains there is no such general tool, but the adaptive
finite element software AFEM of K.~Samuelsson \cite{bsv} has turned out to be effective in a number of cases reported in \cite{bsv}. For the purposes of this paper the so called $hp$-FEM method implemented by H.~Hakula, and first
reported in this paper, is used in several examples with curvilinear boundaries where the previous methods do not apply. The $hp$-FEM method, applied to the computation of moduli of two ring domains previously considered in \cite{bsv} and reported in Tables \ref{table3} and \ref{table4}, provide a significant improvement over the values reported in \cite{bsv}.

For experimental error estimate we have used so called reciprocal identity, which we have not seen used anywhere for the purpose of error estimation. It is our belief that this simple identity should be more widely known. It provides a criterion  for estimating the error of numerical computation of the modulus of a quadrilateral for a very large class of simply connected domains, including those with curved piecewise smooth boundaries. It seems that such a large class of examples has previously not been known for instance in the FEM community. These examples also enable one to experimentally demonstrate the theoretical convergence rates in nontrivial model problem cases as we have shown for the case of $hp$-FEM.

For the very special case of convex polygons with four sides, the modulus of the corresponding quadrilateral is known as an analytic-numeric formula (\ref{hvvmodu}) by \cite{hvv} and this is our starting point. We compare the performance of SC Toolbox, AFEM, and $hp$-FEM against the formula \cite{hvv} and the reciprocal identity test. Next, again using SC Toolbox, AFEM, and $hp$-FEM, we consider polygonal quadrilaterals with more sides, L-shaped quadrilaterals and carry out similar comparision, using again the reciprocal identity as the test quantity. Thereafter, we discuss, now using AFEM, and $hp$-FEM, two classical cases of ring domains, the square frame and the cross in square ring domains previously considered
e.g. in \cite{bsv} where further references may be found. The error estimate in the square frame case is based on the well-known formula (\ref{sqframe}) whereas for the cross in square case we use SC Toolbox and the results from \cite{bsv} as the comparision data. Finally, we also consider several cases of quadrilaterals with curvilinear boundaries, now only using the $hp$-FEM method, because the other methods mentioned above do not apply.

{\bf Acknowledgments.}
We are indebted to the referees for very valuable sets of suggestions concerning the presentation of the results and the use of the SC Toolbox. We also thank Prof. N.~Papamichael for his helpful comments on this paper.

\end{document}